\documentclass[onefignum,onetabnum]{siamonline250211}
\usepackage{amsfonts, amsmath, amssymb}
\usepackage{pgfplots}
\usepackage{theoremref}
\usepackage{mathtools}
\usepackage{commath}
\usepackage{microtype}

\usepgfplotslibrary{groupplots}
\usepgfplotslibrary{fillbetween}

\newtheorem{defn}[theorem]{Definition}
\newtheorem{remark}[theorem]{Remark}
\newcommand{\R}[1]{\mathbb{#1}}
\newcommand{\cl}[0]{\mathcal{L}_{ \mathsf{CV}}^f} 

\newcommand{\ml}{\mathcal{L}_\mathsf{ML}^f} 
\newcommand{\T}{\intercal} 
\newcommand{\mle}[1]{\hat{#1}_{\mathsf{ML}}^f} 
\newcommand{\cle}[1]{\hat{#1}_{\mathsf{CV}}^f} 

\DeclareMathOperator*{\argmin}{arg\,min}
\DeclareMathOperator{\supp}{supp}

 \ifpdf
\hypersetup{
  pdftitle={LENGTHSCALE ESTIMATION FOR THE GAUSSIAN KERNEL},
  pdfauthor={S. Pramanick and T. Karvonen}
}
\fi

\title{Lengthscale Estimation for the Gaussian Kernel and Deterministic Response Functions \thanks{\funding{This work was supported by the Research Council of Finland grants 359183 (“Flagship
of Advanced Mathematics for Sensing, Imaging and Modelling”) and 368086 (“Inference
and approximation under misspecification”), as well as the Finnish Ministry of Education and Culture's Pilot for Doctoral Programmes (Pilot project ``Mathematics of Sensing, Imaging and Modelling'').
TK acknowledges the research environment provided by ELLIS Institute Finland.}}}
\author{Subhendu Pramanick\thanks{School of Engineering Sciences, Lappeenranta--Lahti University of Technology LUT,
Finland (\email{subhendu.pramanick@lut.fi}).}
\and Toni Karvonen\thanks{School of Engineering Sciences, Lappeenranta--Lahti University of Technology LUT,
Finland (\email{toni.karvonen@lut.fi}).}}
\headers{Lengthscale Estimation for the Gaussian Kernel}{Subhendu Pramanick and Toni Karvonen}

\begin{document}

\maketitle

\begin{abstract}
The uncertainty estimates Gaussian processes output rely on the chosen covariance kernel and its hyperparameters.
We study the fixed-domain asymptotics of the maximum likelihood and cross-validation estimators of the lengthscale parameter of the popular infinitely smooth Gaussian kernel when the observations are produced by a fixed univariate response function.
No theoretical results for Gaussian lengthscale estimation exist in this setting.
We prove lower bounds on the lengthscale estimators that depend on the regularity of the response function, as measured by the rate of decay of its Fourier transform.
If the response function lives in the Hilbert space of the Gaussian kernel with lengthscale $\ell$, the lower bounds have a straightforward dependency on $\ell$.
For bandlimited functions and functions whose Fourier transforms decay as $\exp(-a \lvert \omega \rvert^\gamma)$ for $a > 0$ and $\gamma > 2$, we prove that the estimators blow up essentially at least as fast as $(\log n)^{(\gamma-2)/(2\gamma)}$ in the number of data points, $n$.
Bandlimited functions are identified with $\gamma = \infty$.
Unlike analogous results for finitely smooth Matérn kernels, our results make no assumptions on the geometry of the observation points.
Numerical experiments illustrate the results.
\end{abstract}



\section{Introduction}
 
Bayesian nonparametric models are widely used in statistics~\cite{castillo2024bayesiannonparametricstatisticsstflour}, machine learning \cite{10.1145/3291044}, as well as engineering and science \cite{kvam2022nonparametric}.
Gaussian process (GP) models are one of the most popular classes of such nonparametric models, in large part due to their computational tractability and their capabilities to output straightforward uncertainty quantification.
A stochastic process is a GP if its finite-dimensional marginals follow the multivariate normal distribution that is uniquely defined by a mean function and a covariance kernel.
GPs can adapt to complex, non-linear patterns in data.
The mean and covariance are usually equipped with hyperparameters that need to be set arbitrarily or estimated from the data. 
Careful selection of the covariance and its hyperparameters is particularly important because they control the magnitude of the posterior variance that quantifies predictive uncertainty.
This article studies the maximum likelihood (ML) and cross-validation (CV) estimators of the lengthscale parameter $\ell > 0$ of the Gaussian covariance kernel, which is given by
\begin{equation} \label{gaussian_kernel}
    K_\ell(x, y) = \exp\bigg( - \frac{ \lVert x - y \rVert^2}{2\ell^2} \bigg) 
\end{equation}
for $x, y \in \mathbb{R}^d$.
Despite its drawbacks~\cite{KarvonenOates2026}, the Gaussian kernel is one of the most popular kernels in GP modelling.

 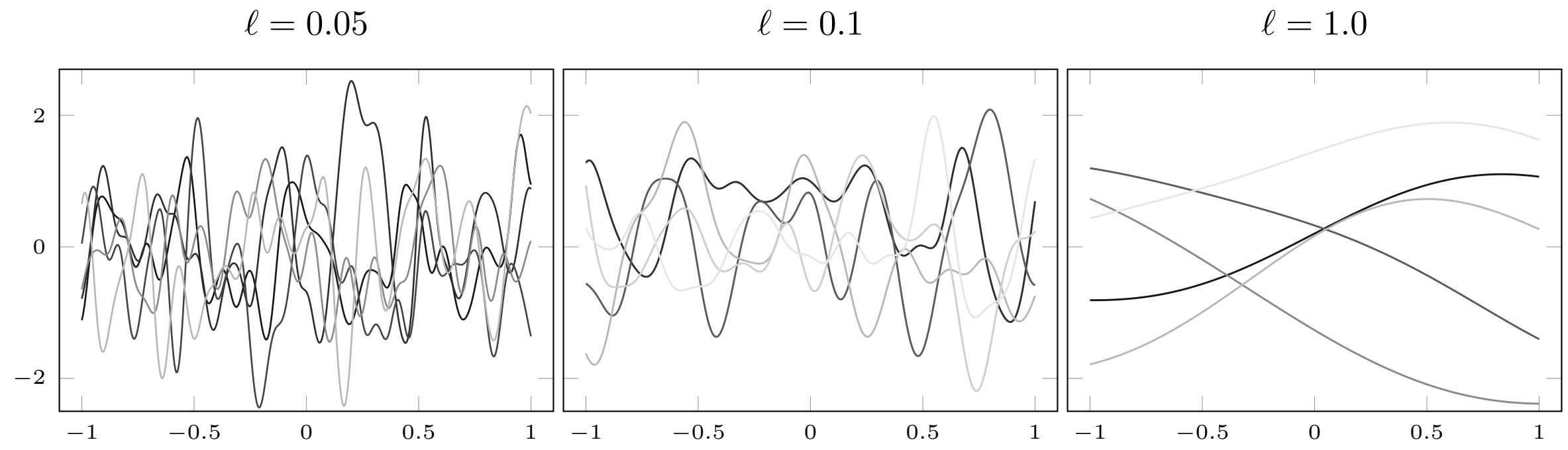
\begin{figure}[t]
    \centering
    \begin{tikzpicture}
        \begin{groupplot}[
            group style={
                group size=3 by 1, 
                horizontal sep=0.1cm,
                y descriptions at= edge left
            },
            width=0.417\linewidth, 
            height=0.32\linewidth,
            xmin=-1.1, xmax=1.1,
            ymin= -2.5, ymax=2.7,
            tick label style={font={\tiny \bfseries}},
            label style={font={\small \bfseries}},
            title style={font={\small \bfseries}},
            legend style={font=\tiny},
            no markers,
        ]
            
            \nextgroupplot[
                title={$\ell = 0.05$}
            ]
            \foreach \n in {1,2,3,6,8}{
            \pgfmathtruncatemacro{\colweight}{9*\n} 
    
            \edef\mycolor{white!\colweight!black}
    
            \edef\temp{
            \noexpand\addplot+[color=\mycolor, line width=0.50pt] 
            table [x=x, y=y\n]{Data/gp_samples0.05.dat};
            }
            \temp
            }

            \nextgroupplot[
                title={$\ell = 0.1$}
            ]
            \foreach \n in {2,4,8, 9, 10}{
            \pgfmathtruncatemacro{\colweight}{9*\n} 
    
            \edef\mycolor{white!\colweight!black}
    
            \edef\temp{\noexpand
                \addplot+ [color=\mycolor, line width=0.55pt] table [x=x, y=y\n]{Data/gp_samples0.1.dat};}
                \temp
            }
            \nextgroupplot[
                title={$\ell = 1.0$}
            ]
            \foreach \n in {1, 4, 6, 8, 10}{
            \pgfmathtruncatemacro{\colweight}{9*\n} 
    
            \edef\mycolor{white!\colweight!black}
            \edef\temp{\noexpand
                \addplot+ [color=\mycolor, line width=0.55pt] table [x=x, y=y\n]{Data/gp_samples1.0.dat};}
                \temp
            }
        \end{groupplot}
    \end{tikzpicture}
    \caption{Sample paths from zero-mean Gaussian processes with Gaussian covariance kernel that has lengthscale $\ell \in \{0.05, 0.1, 1\}$. A small $\ell$ yields highly oscillating sample paths while a large $\ell$ produces slowly varying functions.}
    \label{fig:GP_at_different_lengthscale}
\end{figure}
 
Both the Gaussian kernel and the samples from the GP it induces are analytic (i.e., their Taylor series converge). 
Figure~\ref{fig:GP_at_different_lengthscale} shows how the lengthscale affects samples from the prior.
Choosing a ``good'' lengthscale is crucial, as the lengthscale affects the rate of decay of spatial correlations in the GP prior. 
Selecting too small a lengthscale leads to under-smoothing or overfitting where the GP posterior mean quickly reverts to the prior mean, which results in an under-confident posterior variance. 
On the other hand, too large a lengthscale leads to over-smoothing or underfitting, missing localised non-linear trends and generating over-confident uncertainty quantification.
We study the fixed-domain asymptotics of lengthscale estimators in a frequentist setting.
That is, we assume that the data are produced by an unknown but fixed response function $f$ from which observations are obtained in a fixed bounded set.
We show that the behaviour of the ML and CV estimators depends on the regularity of $f$, as measured by the rate of decay of its Fourier transform. 
Because the approximation theoretic results our analysis relies on do not yet exist in higher dimensions, we focus on the univariate case.

\subsection{Related work} Posterior contraction rates for GP priors with the Gaussian kernel have been studied extensively. It has been shown that with an appropriate choice of lengthscale---contingent on the smoothness of the underlying target function---the posterior distribution contracts with a nearly optimal rate.
These contractions have been investigated in non-parametric regression settings \cite{castillo2014thomas, van2007bayesian, van2009adaptive} and sequence space model settings \cite{pati2015adaptive, doi:10.1137/19M1253010}. Hyperparameter estimation is therefore essential.

When the response function is assumed to be a GP, Loh and Lam \cite{loh2000estimating} showed that a sieve ML estimator for the lengthscale is strongly consistent and asymptotically normal under mild lattice conditions on the observation points.
Loh~\cite{loh2026consistent} proved additional consistency results for scale and lengthscale estimators when the observation points are quasi-uniform on $[0, 1]^d$.
Recently, Qaqish and Li~\cite{qaqish2026sharp} showed that the ML estimator of the lengthscale converges with the rate $\smash{(\log n/ \log \log n)^{-(d+2)/2}}$ in probability if the data are noisy.
In their work, they assumed that the observation points in the input space are sufficiently dense across at least one fixed sub-region of the domain. 
The only work on the Gaussian kernel in the frequentist setting that we are aware of is that of Xu and Stein \cite{xu2017maximum}.

They considered a deterministic response function $f(x)=x^p$, proving asymptotic results for the ML estimator of the scale and numerically studying that of the lengthscale.
A GP whose covariance is a Matérn kernel of regularity $\nu > 0$ has samples that are, essentially, $\nu$ times weakly differentiable.
The Gaussian kernel is obtained at the limit $\nu \to \infty$.
The behaviour of various estimators of the regularity and other parameters of the Matérn class is well-understood~\cite{karvonen2023asymptotic, KarvonenBachoc2025, karvonen2020maximum, korte2025smoothness, Loh2005, loh2021fixed, petit2025asymptotic, Ying1991, zhang2004inconsistent}, both when $f$ is a GP and when it is a fixed deterministic function.
In particular, as noted in \cite[Sec. 3.4]{karvonen2023asymptotic}, the ML estimator of $\nu$ tends to infinity if the response function is infinitely smooth.
This is somewhat reminiscent of the results of this article, which we review next.

\subsection{Our contributions} 

We assume that $f \colon \mathbb{R} \to \mathbb{R}$ is a fixed deterministic function whose values are observed at points $X_n = \{x_1, \ldots, x_n\} \subset [-1, 1]$.
We consider the ML and CV estimators $\mle{\ell}(X_n)$ and $\cle{\ell}(X_n)$ of the lengthscale $\ell$ (see Section~\ref{sec:prelim} for definitions).
The following theorem collects our main results.
More detailed versions appear as Theorems~\ref{th:BdSubsetRKHS} to~\ref{th:PowerExpo}.
The reproducing kernel Hilbert space (RKHS) of the Gaussian kernel with lengthscale $\ell$ is denoted $\mathcal{H}_\ell$ and the Fourier transform of $f$ is $\hat{f}$.

\begin{theorem}\label{intro_thm}
Let $X_n \subset [-1, 1]$ be any sets of $n$ pairwise distinct points.
Suppose that the lengthscale estimators are computed by optimising over $[\delta, \infty)$ for some $\delta > 0$.
    \begin{enumerate}
        \item[(i)] If $f \in \mathcal{H}_{{\ell}}$, then
        \[
        \liminf_{n \to \infty} \mle{\ell}(X_n)\geq {\ell} \quad \text{ and } \quad  \liminf_{n \to \infty}\cle{\ell}(X_n)\geq \frac{\ell}{4}.
        \]
        \item[(ii)] If $f$ is bandlimited and satisfies $\lVert \hat{f} \rVert_{L^2(\mathbb{R})} < \infty$, then there are $C_1, C_2 > 0$ such that, for all sufficiently large $n \in \mathbb{N}$,
        \[
        \mle{\ell}(X_n) \geq C_1 (\log n)^{\frac{1}{2}} \quad \text{ and } \quad \cle{\ell}(X_n) \geq C_2 (\log n)^{\frac{1}{2}} .
        \]
        \item[(iii)] Let $\gamma>2$. 
        If $\lvert \hat{f}(\omega) \rvert \leq M\exp(-a|\omega|^\gamma)$ for some $a, M > 0$, then there are $C_3, C_4 > 0$ such that, for all sufficiently large $n \in \mathbb{N}$,
        \begin{align*}
        \mle{\ell}(X_n) 
        \geq C_3 (\log n)^{\frac{\gamma-2}{2\gamma}} 
        \quad \text{ and } \quad \cle{\ell}(X_n) 
        \geq C_4 (\log n)^{\frac{\gamma-2}{2\gamma}} .
        \end{align*}
    \end{enumerate}
\end{theorem}

Result~(i) is a weaker deterministic version of the results in~\cite{loh2026consistent, loh2000estimating, qaqish2026sharp} on the consistency of lengthscale estimators.
The Gaussian RKHSs are nested: $\mathcal{H}_\eta \subset \mathcal{H}_\ell$ if $\eta \geq \ell$.
This means that, for the ML estimator, Result~(i) is the best one could hope for as $f \in \mathcal{H}_\ell$ does not exclude $f \in \mathcal{H}_\eta$ for $\eta > \ell$.
Nothing resembling Results~(ii) and~(iii) exists in the literature to the best of our knowledge.
A function $f$ is in $\mathcal{H}_\ell$ if $\int_\mathbb{R} \lvert \hat{f}(\omega) \rvert^2 \exp(\ell^2 \omega^2 / 2) \dif \omega < \infty$ (see Section~\ref{sec:rkhs}).
Results~(ii) and~(iii) apply to functions that are smoother than the Gaussian kernel, being in $\mathcal{H}_\ell$ for every $\ell > 0$.
The results state that for such functions the lengthscale estimators blow up with rates that may depend on the behaviour of the Fourier transform.
Because a bandlimited function with bounded Fourier transform satisfies $\hat{f}(\omega) = \mathcal{O}(\exp(-a|\omega|^\gamma))$ for every $\gamma > 0$, Result~(ii) can be thought of as a special case of Result~(iii) with $\gamma = \infty$.
Note that we assume nothing about the observation points $X_n$. 
In particular, the points need not be space-filling.
This is essentially due to the identity theorem for real-analytic functions, which states that a univariate real-analytic function is uniquely determined by its values on a sequence of points that has an accumulation point~\cite[Sec.\@~1.2]{KrantzParks2002}. 
In contrast, all parameter estimation results for finitely smooth Matérns require that the observation points be space-filling if $f$ is deterministic~\cite{karvonen2023asymptotic, KarvonenBachoc2025, karvonen2020maximum, korte2025smoothness}.

The main ingredients in our proofs are estimates on the log determinant of the Gram matrix for the Gaussian kernel. 
The results for the ML estimator use recent bounds by Floater~\cite{floater2026bounds}.
Those for the CV estimator are based on a connection between the determinant of the Gram matrix and the posterior variance.
This makes it possible to leverage upper and lower bounds on the variance from~\cite{karvonen2025approximation}.

\section{Setting and preliminaries} \label{sec:prelim}

This section defines the setting we consider and reviews some auxiliary results.
In this section we consider Gaussian processes on an arbitrary infinite set $\Omega$, called the \textit{domain}.
A kernel $K \colon \Omega\times\Omega \to \mathbb{R} $ is \emph{positive-semidefinite} if the Gram matrix  $K(X_n)=(K(x_i,x_j))_{i,j=1}^n$ is positive-semidefinite for any $n \in \R{N}$ and any set of points $X_n=\{x_1, x_2,\dots, x_n\} \subset \Omega$.
That is,
\[
\alpha^{\intercal}K(X_n)\alpha =\sum_{i=1}^{n}\sum_{j=1}^{n} \alpha_{i}\alpha_{j}K(x_i, x_j)\geq0 
\]
for any $n \in \mathbb{N}$, any $x_1, \ldots, x_n \in \Omega$, and any $\mathbf{\alpha}=(\alpha_1, \alpha_2, \dots, \alpha_n)\in \R{R}^n$.
A kernel is positive-definite if the Gram matrix is \emph{positive-definite} whenever $x_1, \ldots, x_n$ are pairwise distinct.

\subsection{Gaussian processes}

Let $m \colon \Omega \to \R{R}$ be a function and $K_{\theta} \colon \Omega \times \Omega \to \R{R}$ a positive-definite kernel parameterised by parameters $\theta$ that lie in a set $\Theta$ of feasible parameters.
A Gaussian process $\mathsf{f}_\theta$ with mean $m$ and covariance $K_\theta$ is a stochastic process such that for any points $X_n = \{x_1, \ldots, x_n \}$ the random vector
\[
\big(\mathsf{f}_\theta(x_1), \ldots, \mathsf{f}_\theta(x_n)\big)
\]
is an $n$-dimensional normal random variable with mean $m(X_n) = (m(x_i))_{i=1}^n \in \R{R}^n$ and covariance matrix $K_\theta(X_n) \in \R{R}^{n \times n}$. 
We write $\mathsf{f}_\theta \sim \mathcal{GP}(m, K_{\theta})$ because a Gaussian process is uniquely defined by the mean and covariance functions, which satisfy
\begin{equation*}
    m(x) = \mathbb{E}[\mathsf{f}_\theta(x)] \quad \text{ and } \quad K_\theta(x, y) = \mathrm{Cov}(\mathsf{f}_\theta(x), \mathsf{f}_\theta(y)) 
\end{equation*}
for all $x, y \in \Omega$.
The covariance determines properties such as smoothness and periodicity of a Gaussian process.
For example, the samples of a zero-mean Gaussian process with the Matérn covariance 
\begin{equation*}
    K(x, y) = \frac{2^{1-\nu}}{\Gamma(\nu)} \bigg( \frac{\sqrt{2\nu} \lVert x - y \rVert}{\ell} \bigg)^\nu \mathcal{K}_\nu\bigg( \frac{\sqrt{2\nu} \lVert x - y \rVert}{\ell} \bigg)
\end{equation*}
with smoothness $\nu > 0$ on $\Omega = \R{R}^d$ are essentially $\nu$ times weakly differentiable.
Here $\mathcal{K}_\nu$ is the modified Bessel function of the second kind.
The Gaussian kernel~\eqref{gaussian_kernel}, which is obtained from Matérns at the limit $\nu \to \infty$, defines a Gaussian process whose samples are infinitely differentiable.
For simplicity, from now on we tacitly assume that all Gaussian processes have zero mean (i.e., $m = 0$).
See~\cite{gine2021mathematical,williams2006gaussian} or \cite[Chapter 5]{kanagawa2025gaussian} for introductions to Gaussian processes.

\begin{figure}[htbp]
    \centering
    \begin{tikzpicture}
        \begin{axis}[
            width=1.0\linewidth, 
            height=0.40\linewidth,
            xmin=-0.10, xmax=5.1,
            ymin=-0.60, ymax=1.8,
            grid=none,
            tick label style={font={\tiny \bfseries}},
            label style={font={\small \bfseries}},
            title style={font={\small \bfseries}},
            legend style={font=\tiny},
        ]
    \addplot[name path=upper, draw=none, forget plot] table{Data/Var_Upper_limit.dat};
    \addplot[name path=lower, draw=none, forget plot] table{Data/Var_Lower_limit.dat};
    \addplot [fill= red!20, opacity=0.8] fill between [of=upper and lower];
    \addplot [line width = 0.6pt, color= cyan] table{Data/Mean_prediction.dat};
    \addplot[color=black, only marks, mark = *, mark size = 1pt, line width=0.02pt, forget plot] table{Data/True_response7.dat};
     \addplot[color=orange, line width=0.6pt] table{Data/True_response280.dat};
     \legend{$95\%~$Credible interval, Posterior mean, True response function}
    \end{axis}
\end{tikzpicture}
\caption{\emph{Gaussian process interpolation.} The Gaussian process posterior mean (cyan line) given in~\eqref{eq:gp-mean} is computed based on observed data (black dots) from a true response function (orange line). The shaded light-red area represents the 95\% credible interval, calculated directly from the posterior variance~\eqref{eq:gp-variance} of the Gaussian process.}
\label{fig:UQ}
\end{figure}
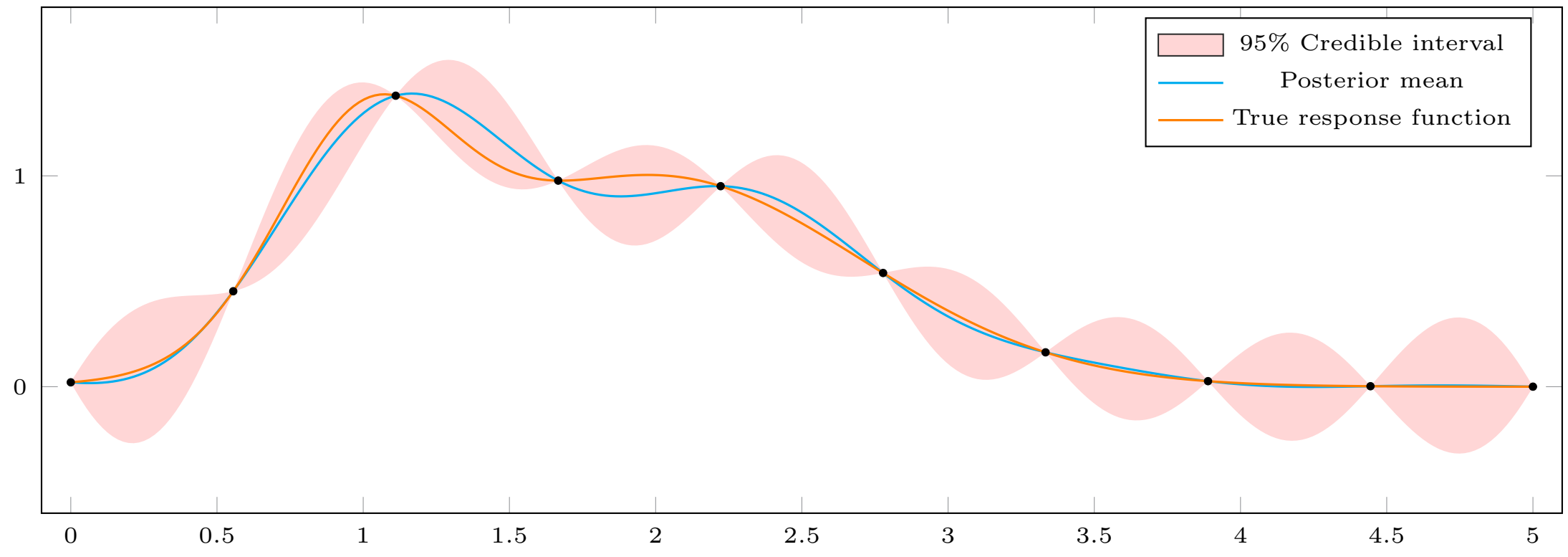
Suppose that a deterministic response function $f \colon \Omega \to \R{R}$
is modelled as a Gaussian process $\mathsf{f}_\theta \sim \mathcal{GP}(0, K_{\theta}) $. Conditioning this process on the exact evaluations (i.e., data)
\[
f(X_n) = \big(f(x_1), \ldots, f(x_n)\big) \in \R{R}^n
\] 
of $f$ at some pairwise distinct points $X_n$ yields a conditional Gaussian process with mean
\begin{equation} \label{eq:gp-mean}
\mu_{\theta,f}(x \mid X_n) 
= \R{E}[\mathsf{f}_\theta(x) \mid X_n, f(X_n)] 
= K_{\theta}(x, X_n)^\T K_{\theta}(X_n)^{-1} f(X_n)
\end{equation}
and variance
\begin{equation} \label{eq:gp-variance}
V_{\theta}(x \mid X_n) 
= \mathrm{Var}[\mathsf{f}_\theta(x) \mid X_n, f(X_n)] 
= K_{\theta}(x, x) - K_{\theta}(x, X_n)^\T K_{\theta}(X_n)^{-1} K_{\theta}(x, X_n),
\end{equation}
where $K_{\theta}(x, X_n) = (K_\theta(x, x_i))_{i=1}^n \in \R{R}^n$.
Figure~\ref{fig:UQ} illustrates the conditional Gaussian process.
The conditional mean provides predictions for $f(x)$ at $x \notin X_n$ while the variance quantifies the associated uncertainty, for example via 95\% credible intervals as in Figure~\ref{fig:UQ}.
The mean \emph{interpolates} the data, in that $\mu_{\theta,f}(x_i \mid X_n) = f(x_i)$ for every $i \in \{1, \ldots, n\}$.
However, observe that as such the variance in~\eqref{eq:gp-variance} does not depend on $f$.
If the Gaussian process $\mathsf{f}_\theta$ is not a good model for $f$, uncertainty quantification may therefore easily become unreliable, meaning that the conditional variance is not an accurate surrogate of the true interpolation error~\cite{karvonen2020maximum}.
This can be combated by estimating the covariance parameters $\theta$ from the data.

From the perspective of Bayesian hierarchical modelling, Gaussian process interpolation is a two-level model.
Let $\mathcal{D}_n = (X_n, f(X_n))$ denote the data.
At the lower level, the prior $p( \mathsf{f}_\theta)$ over the latent response function $f$ is governed by a Gaussian process $\mathsf{f}_\theta\sim \mathcal{GP}(0, K_\theta)$ with parameters $\theta \in \Theta$.
Through Bayes' theorem, conditioning on the data yields the posterior
\[
 p( \mathsf{f}_\theta \mid \mathcal{D}_n, \theta) = \frac{p( \mathcal{D}_n \mid \mathsf{f}_\theta, \theta) p({\mathsf{f}_\theta} \mid \theta)} {p(\mathcal{D}_n \mid \theta)},
\] 
where
\begin{equation*}
    p(\mathcal{D}_n \mid \theta) = \int p( \mathcal{D}_n \mid \mathsf{f}_\theta, \theta) p({\mathsf{f}_\theta} \mid \theta) \dif{\mathsf{f}_\theta}
\end{equation*}
is the marginal likelihood.\footnote{Because we assume that the observations are not corrupted by noise, the likelihood $p(\mathcal{D}_n \mid \mathsf{f}_\theta)$ is a Dirac distribution and lacks a density function, which precludes the use of Bayes' theorem~\cite[Section~2.5]{Cockayne2019}. However, the posterior can be derived from the usual Gaussian conditioning results and coincides with the expression that one gets by setting the noise variance to zero in the posterior obtained from Bayes' theorem under Gaussian observation noise.}
At the second level are the covariance parameters $\theta$ that are hyperparameters for the Gaussian process model.
A prior could be placed over the hyperparameters, but this is often both computationally and theoretically intractable. 
Selecting hyperparameters that maximise the marginal likelihood is more tractable.
This approach is known as type II maximum likelihood~\cite[Chapter~5]{williams2006gaussian}.
Because the prior for $f$ is a Gaussian process, it is straightforward to work out that
\begin{equation} \label{eq:marginal-likelihood}
    \log p(\mathcal{D}_n \mid \theta) = -\frac{1}{2} f(X_n)^\T K_\theta(X_n)^{-1} f(X_n) - \frac{1}{2} \log \det K_\theta(X_n) - \frac{n}{2} \log 2\pi .
\end{equation}
The maximum likelihood estimator of $\theta$ is any maximiser of~\eqref{eq:marginal-likelihood}. 
For our purposes it is convenient to work with a minimiser of the negative log-likelihood, an equivalent quantity.

\begin{defn}[{Maximum likelihood estimator}]\thlabel{mle_def} 

The \emph{maximum likelihood (ML) estimator} $\mle{\theta}(X_n)$ of $\theta$ is any maximiser of~\eqref{eq:marginal-likelihood}.
That is,
\[
\mle{\theta}(X_n) \in \argmin_{\theta \in \Theta } \ml(\theta \mid X_n),
\]
where
\begin{equation} \label{eq:L-ML}
    \ml(\theta \mid X_n) = f(X_n)^\intercal K_{\theta}(X_n)^{-1} f(X_n) + \log \det K_{\theta}(X_n) .    
\end{equation}

We set $\mle{\ell}(X_n) = \infty$ if no minimiser exists.
\end{defn}
For terminological simplicity, we refer to this modified negative log-likelihood, denoted as $\ml(\theta \mid X_n)$, as ML likelihood.

The cross-validation estimator (also known as Geisser's surrogate predictive probability estimator) for Gaussian processes originates in predictive sample reuse methodology~\cite{geisser1975predictive, sundararajan1999predictive}. Specifically, it is used to estimate a parameter by maximising the predictive likelihood of left-out observations.
In a leave-one-out framework, for each $i \in \{1, \dots, n\}$, the $i^{\text{th}}$ sample is omitted to form a reduced set $X_n^{i} = X_n \setminus \{x_i\}$. The cross-validation estimator for $\theta$ maximises the product of the predictive distributions of the omitted samples, $\prod_{i=1}^{n} p(f(x_i) \mid  \mathcal{D}_{n}^i, \theta)$, where $\mathcal{D}_n^i$ is defined analogously to $X_n^i$. 
Due to Gaussianity this product has a simple expression. 
Again, it is more convenient to work with an equivalent minimisation problem.

\begin{defn}[Cross-validation estimator]\label{cv_def}
    The \emph{leave-one-out cross-validation (CV) estimator} $\cle{\theta}(X_n)$ is any maximiser of $\prod_{i=1}^{n} p(f(x_i) \mid  \mathcal{D}_{n}^i, \theta)$. 
    Equivalently,
    \[
    \cle{\theta}(X_n) \in \argmin_{\theta \in \Theta} \cl(\theta \mid X_n),
    \]
    where
    \begin{equation} \label{eq:L-CV}
        \cl(\theta \mid X_n) = \sum_{i=1}^{n}\left[ \frac{(f(x_i) - \mu_{\theta, f}(x_i \mid X_n^i))^2}{V_{\theta}(x_i \mid X_n^i)} + \log V_{\theta}(x_i \mid X_n^i)\right] .
    \end{equation}
    We set $\cle{\ell}(X_n) = \infty$ if no minimiser exists.
\end{defn}

We call $\cl(\theta \mid X_n)$ the CV likelihood.
In addition to being fundamentally linked by their shared likelihood-based origins, these two estimators are also theoretically related via lemmas to be presented in Section~\ref{sec:lemmas} and the theory in~\cite{FongHolmes2020}.
Specifically, Lemma~\ref{lem:DetToVar} establishes a relationship between the posterior variance in cross-validation and the determinant of the Gram matrix in maximum likelihood estimation, the second terms in~\eqref{eq:L-ML} and~\eqref{eq:L-CV}. Furthermore, Lemma~\ref{lem:NormBd} reconciles the data-fitting components, the first terms in~\eqref{eq:L-ML} and~\eqref{eq:L-CV}, of both frameworks. Given these mathematical correspondences, it is not surprising that the two estimators exhibit analogous properties.
We study the ML and CV estimators of the lengthscale parameter $\ell$ of the Gaussian kernel in~\eqref{gaussian_kernel}. 

\subsection{Reproducing kernel Hilbert spaces} \label{sec:rkhs}

Our proofs rely on \emph{reproducing kernel Hilbert spaces} (RKHSs).
An RKHS $\mathcal{H}$ on $\Omega$ over $\R{R}$ is a Hilbert space containing real valued functions on $\Omega$ such that the point evaluation functional
\[  
 h \mapsto h(x)
\] 
is bounded for every $x \in \Omega$.
Each RKHS has a unique positive-semidefinite reproducing kernel $K \colon \Omega \times \Omega \to \R{R}$ that satisfies the following reproducing property for every $h \in \mathcal{H}$ and $x \in \Omega$:
\begin{equation} \label{eq:repro}
h(x) = \langle h , k_x \rangle, \quad \text{ where} \quad k_x(y) \coloneqq K(y,x).
\end{equation}
Conversely, each positive-semidefinite kernel $K$ induces a unique real RKHS on $\Omega$.
This RKHS is denoted by $\mathcal{H}(K)$. 
See~\cite{paulsen2016introduction, berlinet2011reproducing} for an introduction to the theory of RKHSs. 
When the kernel is parameterised by $\theta$, we write $\mathcal{H}_\theta = \mathcal{H}(K_\theta)$ and $\langle \cdot, \cdot \rangle_\theta = \langle \cdot, \cdot \rangle_{\mathcal{H}(K_\theta)}$ for brevity.

The Fourier transform of $f \in L^2(\mathbb{R})$ is
\[
\hat{f}(\omega) = \frac{1}{\sqrt{2\pi}}\int_{\mathbb{R}} f(x)e^{- \sqrt{-1} \omega x} \, \dif x.
\]
The Gaussian kernel that we consider is \emph{translation-invariant} (or stationary), meaning that it can be written as $K(x,y) = \Phi(x-y)$ for any $ x, y \in \R{R}$ and some continuous and integrable function $\Phi$ on $\R{R}$. 
The RKHS of any such translation-invariant kernel admits a convenient spectral representation \cite[Theorem 10.12]{wendland2004scattered}.
Namely,
\[
\mathcal{H}(K, \R{R}) = \left\{ f \in L^2(\R{R}) \cap C(\R{R}) : \int_{\R{R}} |\hat{f}(\omega)|^2|{\hat{\Phi}(\omega)}|^{-1} \dif \omega < \infty \right\}
\]
and the inner product is
\[
{\langle f, g \rangle}_{\mathcal{H}(K, \R{R})} = \frac{1}{\sqrt{2\pi}}  \int_{\R{R}} \frac{\hat{f}(\omega) \overline{\hat{g}(\omega)}}{\hat{\Phi}(\omega)} \dif \omega.
\]
Here the notation $\mathcal{H}(K, \R{R})$ makes it explicit that this is an RKHS on $\R{R}$.
For the Gaussian kernel $K_\ell$ with lengthscale $\ell$ we have
\begin{equation*}
    \Phi(r) = \exp\bigg( - \frac{r^2}{2\ell^2} \bigg), \quad \text{ so that } \quad \hat{\Phi}(\omega) = \ell \exp{\bigg(-\frac{\omega^2 \ell^2}{2}\bigg)} .
\end{equation*}
Therefore the RKHS of the Gaussian kernel with lengthscale $\ell$ is (e.g.,\@~\cite[Theorem 2.4]{kanagawa2025gaussian})
\begin{equation} \label{eq:gaussian-rkhs}
\mathcal{H}(K_\ell, \R{R}) = \left\{ f \in L^2(\mathbb{R}) \cap C(\mathbb{R}) : \|f\|^2_{\ell, \R{R}} = \frac{1}{\sqrt{2\pi}\ell} \int_{\mathbb{R}} |\hat{f}(\omega)|^2 \exp\bigg( \frac{\ell^2 \omega^2}{2} \bigg) \dif\omega < \infty \right\} .
\end{equation}
This RKHS is equipped with the inner product
\begin{equation}\label{eq:gaussian-inner-product}
{\langle f, g \rangle}_{\mathcal{H}(K_\ell, \R{R})} = \frac{1}{\sqrt{2\pi}\ell} \int_{\mathbb{R}} \hat{f}(\omega)\overline{\hat{g}(\omega)}\exp\bigg( \frac{\ell^2 \omega^2}{2} \bigg)\dif\omega .
\end{equation}
We consider the Gaussian kernel on $\Omega = [-1, 1]$, a subset of $\R{R}$ with non-empty interior.
By~\cite[Theorem 5.7]{paulsen2016introduction}, the RKHS of Gaussian kernel on $\Omega$ is the restriction of $\mathcal{H}(K_\ell, \R{R})$ onto $\Omega$. 
Precisely,
\[
\mathcal{H}(K_\ell, \Omega) = \big\{ f_{|\Omega} \colon f \in \mathcal{H}(K_\ell, \R{R}) \big\}
\] 
with norm $\|f\|_{\ell, \Omega}= \inf\{\|g\|_{\ell,\R{R}} \colon f = g_{|\Omega}\}$. 
Interestingly, every $f \in \mathcal{H}(K_\ell, \Omega)$ has a unique extension in $\mathcal{H}(K_\ell, \R{R})$ \cite[Corollary 3.8]{steinwart2006explicit} and the extension is isometric, meaning that $\|f\|_{\ell, \Omega} = \|f\|_{\ell, \R{R}}$. Therefore we can identify each $f \in \mathcal{H}(K_\ell, \Omega)$ with its unique extension in $\mathcal{H}(K_\ell, \R{R})$.
Consequently, in the remainder of the article we do not make distinction between an RKHS on $\R{R}$ and $\Omega$ and simply write 
\begin{equation*}
    \mathcal{H}_\ell = \mathcal{H}(K_\ell, \Omega) = \mathcal{H}(K_\ell, \R{R}) \quad \text{ and } \quad \|f\|_\ell = \|f\|_{\ell,\Omega} = \|f\|_{\ell,\R{R}},   
\end{equation*}
where the Gaussian RKHS and its norm are defined in~\eqref{eq:gaussian-rkhs} and~\eqref{eq:gaussian-inner-product}.
For additional characterisations and properties of RKHS associated with the Gaussian kernel, we refer to \cite{sun2008reproducing, minh2010some, steinwart2006explicit}.


\subsection{Bounds and equivalences} \label{sec:lemmas}

Our main results in Section \ref{sec:MainResults} use bounds on posterior variances and the RKHS norm of the posterior mean. 
Here we review auxiliary results that establish an identity between the Gaussian process posterior variance and the Gram determinant (i.e., determinant of the Gram matrix) and provide bounds for these quantities and the RKHS norm of the posterior mean. 
Results that are not specific to the Gaussian kernel are parameterised by a general kernel hyperparameter~$\theta$.
The results specific to the Gaussian kernel are parameterised by the lengthscale parameter~$\ell$. 

  A well-known lemma (e.g., Lemma~2.1 in~\cite{karvonen2023asymptotic}) connects the determinant of the Gram matrix to the posterior variance. We include an elementary proof here for the sake of completeness.

\begin{lemma}\label{lem:DetToVar}
Define $V_{\theta}(x_1 \mid X_{0}) = K_{\theta}(x_1,x_1)$. Then
    \[
    \det K_{\theta}(X_n) =\prod_{i=1}^{n}  V_{\theta}(x_i \mid X_{i-1}). 
    \] 
\end{lemma}
\begin{proof}
We proceed by recursively expanding the determinant of the covariance matrix.  
For $i \geq 2$, write the covariance matrix of the first $i$ points in block form
\[
K_\theta(X_i)
=
\begin{pmatrix}
K_\theta(X_{i-1}) & K_\theta(X_{i-1},x_i) \\
K_\theta(x_i,X_{i-1}) & K_\theta(x_i,x_i)
\end{pmatrix}.
\]
Let $a = K_\theta(x_i,x_i)$, $b = K_\theta(X_{i-1},x_i)$ and $C = K_\theta(X_{i-1})$.  
Applying the block determinant identity
\[
\det\begin{pmatrix}
C & b \\[2pt]
b^\intercal & a
\end{pmatrix}
= \det(C)\,\bigl(a - b^\T C^{-1} b \bigr)
\]
gives
\begin{align*}
\det K_\theta(X_i) &= \det K_\theta(X_{i-1})\,\big(K_\theta(x_i,x_i)- K_\theta(x_i,X_{i-1})^\T K_\theta(X_{i-1})^{-1}\,K_\theta(X_{i-1},x_i)\big) \\
& = \det K_\theta(X_{i-1})\,V_\theta(x_i \mid X_{i-1}).
\end{align*}
We obtain the claim by repeating the above procedure:
\[
\det K_{\theta}(X_n) =\prod_{i=1}^{n}  V_{\theta}(x_i \mid X_{i-1}) .
\]
\end{proof}
From Lemma~\ref{lem:DetToVar} we obtain
\[
\det K_{\theta}(X_n) = V_{\theta}(x_n \mid X_{n-1})\prod_{i=1}^{n-1}  V_{\theta}(x_i \mid X_{i-1})= V_{\theta}(x_n \mid X_{n-1}) \det K_{\theta}(X_n^n). 
\]
Moreover, a permutation of indices yields
\begin{equation}\label{eq:VarToDet}
    V_\theta(x_i \mid X_n^i) = K_\theta(x_i,x_i) - K_\theta(x_i,X_n^i)^\intercal K_\theta(X_n^i)^{-1}K_\theta(x_i,X_n^i) = \frac{\det K_\theta(X_n)}{\det K_\theta(X_n^i)}.
\end{equation}
Hence the posterior variance at a left-out observation is given by the ratio of the determinant of the Gram matrix to the determinant of its corresponding leave-one-out Gram matrix.

The lemma below expresses the data-fit term in $\ml(\theta \mid X_n)$ as the RKHS norm of the conditional mean.
\begin{lemma}[Lemma~2.2 in \cite{karvonen2023asymptotic}]\label{lem:NormBd}
We have
\[
f(X_{n})^\T K_{\theta}(X_{n})^{-1} f(X_{n})
=\sum_{i=1}^{n} \frac{(f(x_i) - \mu_{\theta, f}(x_i \mid X_{i-1}))^2}{V_{\theta}(x_i \mid X_{i-1})}= \big\| \mu_{\theta, f}(\,\cdot \mid X_{n}) \big\|_{\theta}^{2}.
\]
Moreover, if $f \in \mathcal{H}_\theta$, then
\[
\lVert \mu_{\theta, f}(\,\cdot \mid X_{n}) \rVert_{\theta}^{2}
\leq \| f \|^{2}_{\theta}.
\]
\end{lemma}

Lemma~\ref{lem:NormBd} implies that
    \begin{equation*}
         \frac{(f(x_i) - \mu_{\theta, f}(x_i \mid X_n^i))^2}{V_{\theta}(x_i \mid X_n^i)} = \|\mu_{\theta,f}( \cdot \mid X_n)\|_\theta^2 - \|\mu_{\theta,f}( \cdot \mid X_n^i)\|_\theta^2.
    \end{equation*}
    This yields a simple upper bound on the prediction error term in $\cl(\theta \mid X_n)$ (Definition~\ref{cv_def}):
    \begin{equation}\label{eq:CVNorm}
             \frac{(f(x_i) - \mu_{\theta, f}(x_i \mid X_n^i))^2}{V_{\theta}(x_i \mid X_n^i)} \leq \|f\|_\theta^2,
    \end{equation}
for any $f \in \mathcal{H}_\theta$.
The identities in Lemma~\ref{lem:DetToVar} are based solely on linear algebra.
However, posterior variances can also be understood through RKHS theory. 
In particular, Novak and Woźniakowski~\cite[Section~10.2]{novak2010tractability} characterise the posterior variance as the minimal squared worst-case error in $\mathcal{H}_\theta$ based on function evaluations at $X_n$. Specifically, the posterior variance $V_{\theta}(x \mid X_n) = K_{\theta}(x, x) - K_{\theta}(x, X_n)^\T K_{\theta}(X_n)^{-1} K_{\theta}(x, X_n)$ equals the minimum of the quadratic form
\[
K_\theta(x, x) - 2 \sum_{l=1}^n u_l K_\theta(x, x_l) + \sum_{k,l=1}^n u_k u_l K_\theta(x_k, x_l)
\] over $u_1,\dots,u_n \in \mathbb{R}$. 
With the reproducing property~\eqref{eq:repro}, this optimisation problem is equivalently formulated as
\begin{align*}
    V_\theta(x \mid X_n) &= \inf_{u_1,\dots,u_n \in \mathbb{R}} \bigg[K_\theta(x, x) - 2 \sum_{l=1}^n u_l K_\theta(x, x_l) + \sum_{k,l=1}^n u_k u_l K_\theta(x_k, x_l)\bigg] \\
    & = \inf_{u_1,\dots,u_n \in \mathbb{R}} \sup_{0 \neq f \in \mathcal{H}(K_\theta)}\frac{| f(x) -\sum_{i=1}^n f(x_{i}) u_{i} |^2}{\|f\|_\theta^2}.
\end{align*}

Karvonen and Suzuki~\cite[Section~3]{karvonen2025approximation} have established upper and lower bounds on the worst-case error in the RKHS of the Gaussian kernel. 
Through the worst-case interpretation of the posterior variance we obtain the following theorem.

\begin{theorem}[Theorem 1.1 in~\cite{karvonen2025approximation}]\label{th:WCE_Bound}
Consider the Gaussian kernel in~\eqref{gaussian_kernel}.
If $x_1, \dots, x_n \in [-1,1]$ are any pairwise distinct points, then for all $x \in [-1,1]\setminus \{x_1, \dots, x_n\}$,
\[
  \frac{c_{1}(\ell)^2}{n!}\left(\frac{1}{4\ell}\right)^{2n} \leq \frac{V_{\ell}(x \mid X_n)}{|\prod_{i=1}^{n}(x-x_i)|^2} \leq \frac{c_{2}(\ell)^2}{n!}\,n^{-1/4}e^{\frac{2\sqrt{n}}{\ell}}\bigg(\frac{1}{\ell}\bigg)^{2n},
\] 
where
\begin{align*}
    c_{1}(\ell) &= \sqrt{\frac{2\ell^2}{1 + \ell^2}}\, e^{-\frac{1}{32\ell^4}}e^{-\frac{1}{2\ell^2}}, \\
    c_{2}(\ell) &=\Bigg[ \frac{ \exp\big(\frac{1}{2\ell^2} + \frac{2}{\ell}\big) }{ 2\sqrt{\pi} \big(\frac{1}{\ell}\big)^{1/2} } \bigg(3 + \sqrt{\frac{\pi}{e}} + 2\sqrt{\frac{\ell}{\pi}} + \frac{7}{16}\ell + \frac{1}{\ell}\Big(2 + \frac{1}{\ell^2}\Big) \bigg)\Bigg]^{\frac{1}{2}}.
\end{align*}
\end{theorem}

The Vandermonde polynomial \(V_{X_n}\), the determinant of the Vandermonde matrix \((x_{i}^{j-1})_{i,j=1}^n\), satisfies the identity
\[
V_{X_n}=\prod _{1\le i<j\le n}(x_{j}-x_{i}).
\]
Let $p_n = \prod_{i=1}^{n-1} (i!)$.
Theorem~\ref{th:WCE_Bound} implies two-sided bounds for the Gaussian Gram determinant.

\begin{lemma}\label{lem:DetBdUsngWCE}
    Let $c_1$ and $c_2$ be the constants in Theorem~\ref{th:WCE_Bound}. Then
    \begin{align*}
     \Big(\frac{1}{4\ell}\Big)^{n(n-1)}\frac{c_1(\ell)^{2(n-1)}}{p_n}  \leq \frac{\det(K_\ell(X_n))}{V_{X_n}^2}
     \leq \bigg(\frac{1}{\ell}\bigg)^{n(n-1)}\frac{c_2(\ell)^{2(n-1)}}{p_n((n-1)!)^{1/4}}  \exp\bigg({\frac{2\sum_{i=1}^{n-1}\sqrt{i}}{\ell}}\bigg).
    \end{align*}
\end{lemma}
\begin{proof}
    From Lemma~\ref{lem:DetToVar} we know that $\det K_{\ell}(X_n) = \prod_{i=1}^{n}  V_{\ell}(x_i \mid X_{i-1})$.
    Moreover, note that $V_{X_n} = \prod_{i=1}^{n}\prod_{j< i} (x_i -x_j)$.
    We obtain the claimed bounds by using bounds on $V_\ell(x \mid X_n)$ in Theorem~\ref{th:WCE_Bound}.
\end{proof}

A recent result by Floater~\cite{floater2026bounds} provides a determinant bound for exponential matrices $(e^{x_iy_j})_{i,j=1}^n$, which yields a bound on the Gaussian Gram determinant as a corollary. To present this result we introduce the quantities
\begin{equation}\label{eq:AverageMeanSquareDiff}
\Bar{x} = \frac{1}{n}\sum_{i=1}^{n} x_i \quad \text{ and } \quad S_n = \sum_{i=1}^{n}(x_i - \Bar{x})^2.
\end{equation}
Since the permanent of a matrix $(e^{-(x_i-x_j)^2})_{i,j=1}^n$ ($x_i \in \mathbb{R})$ is bounded above by $n!$,  bounding the permanent of $K_\ell(X_n)$  by $n!$ in~\cite[Corollary~1]{floater2026bounds} yields the lemma below.

\begin{lemma} \label{lem:DetBdFloater}
Consider the Gaussian kernel given in~\eqref{gaussian_kernel}. 
If $x_1, \dots, x_n \in [-1,1]$ are any pairwise distinct points, then 
    \[
    \frac{1}{p_n \ell^{n(n-1)}} e^{- S_n/\ell^2} \leq \frac{\det K_\ell(X_n)}{V_{X_n}^2} \leq  \frac{1}{p_n \ell^{n(n-1)}}.
    \]
\end{lemma}
\begin{remark}
    The two-sided bounds in Lemmas~\ref{lem:DetBdUsngWCE} and~\ref{lem:DetBdFloater} are readily comparable due to their similar structure. Although both have a similar growth rate with respect to $\ell$, the upper bound in Lemma~\ref{lem:DetBdUsngWCE} grows faster by approximately $e^{4n^{3/2}/(3\ell)}$.
\end{remark}

\section{Main results}\label{sec:MainResults}

In this section we prove the results that were summarised as Theorem~\ref{intro_thm} in the introduction.
Each of these results provides lower bounds on the ML and CV estimators of the lengthscale of the Gaussian kernel. 
The proofs are based on a general method for finding lower bounds for hyperparameter estimators in Gaussian process modelling, introduced in~\cite[Theorem 2.3]{karvonen2023asymptotic}. 
Below we provide a version of this method that best suits our purposes.

\begin{lemma}\label{lem:LikelihoodComp}
      Suppose that $K_\theta$ is parameterised by $\theta \in \R{R}_+$. Let $\theta_0, \theta > 0$.  
    \begin{enumerate}
        \item If $f$ is any function in $\mathcal{H}_{\theta_0}$ such that 
    \begin{equation}
       \| f \|^{2}_{\theta_0} + \log \frac{\det K_{\theta_0}(X_n)}{\det K_{\theta}(X_n)} < 0,
    \end{equation} 
     then $\ml(\theta_0 \mid X_n) < \ml(\theta \mid X_n)$.
     
        \item If $f$ is any real valued function on $\Omega$ 
        such that
    \begin{equation}
       \sum_{i=1}^{n} \bigg[ \frac{(f(x_i) - \mu_{\theta_0, f}(x_i \mid X_n^i))^2}{V_{\theta_0}(x_i \mid X_n^i)} + \log \frac{V_{\theta_0}(x_i \mid X_n^i)}{ V_{\theta}(x_i \mid X_n^i)} \bigg] < 0,
    \end{equation} 
     then $\cl(\theta_0 \mid X_n) < \cl(\theta \mid X_n)$.
    \end{enumerate}
\end{lemma}
\begin{proof}
    We begin with
    \begin{align*}
        \ml(\theta_0 \mid X_n) &= f(X_n)^\T K_{\theta_0}^{-1}(X_n) f(X_n) + \log \det K_{\theta_0}(X_n)\notag\\
                    &\leq \| f \|^{2}_{{\theta_0}} + \log \det K_{\theta_0}(X_n)\tag{by Lemma~\ref{lem:NormBd}}\\
                    &\leq \| f \|^{2}_{{\theta_0}} + ||\mu_{{\theta},f}(\cdot \mid X_n)||_{\theta}^2 + \log \det K_{\theta_0}(X_n)\notag\\
                    &=  \ml(\theta \mid X_n) + \| f \|^{2}_{{\theta_0}}  +\log \frac{\det K_{\theta_0}(X_n)}{\det K_{\theta}(X_n)}.
    \end{align*} 
    Therefore, if 
    \begin{equation}
        \bigg[ \| f \|^{2}_{\theta_0} + \log \frac{\det K_{\theta_0}(X_n)}{\det K_{\theta}(X_n)} \bigg] < 0,
    \end{equation} 
    then $\ml(\theta_0 \mid X_n) < \ml(\theta \mid X_n)$. The proof is analogous in the case of CV likelihood.
\end{proof}
\begin{theorem}\label{th:Gen_Method}
    Suppose that $K_\theta$ is parameterised by $\theta \in \R{R}_+$. Let $0<\theta_1 \leq \theta_0$.  
    \begin{enumerate}
        \item If $\mathcal{B}$ is a subset of $\mathcal{H}_{\theta_0}$ such that
    \begin{equation}\label{gen_cond}
        \limsup_{n\to \infty} \sup_{\theta<\theta_1}\sup_{f \in \mathcal{B}}\bigg[ \| f \|^{2}_{\theta_0} + \log \frac{\det K_{\theta_0}(X_n)}{\det K_{\theta}(X_n)} \bigg] < 0,
    \end{equation} 
     then there exists $n_0 \in \R{N}$ such that $\ml(\theta_0 \mid X_n) < \ml(\theta \mid X_n)$ for every $\theta \in (0, \theta_1)$ and $f \in \mathcal{B}$ for all $n \geq n_0$.
     Hence, $\inf_{f \in \mathcal{B}}\mle{\theta}(X_n) \geq \theta_1 $ when $n \geq n_0$.
        \item If $\mathcal{B}$ is a subset of $\mathcal{H}_{\theta_0}$ such that
    \begin{equation}\label{gen_cond2}
        \limsup_{n\to \infty} \sup_{\theta<\theta_1}\sup_{f \in \mathcal{B}}\sum_{i=1}^{n}\bigg[ \frac{(f(x_i) - \mu_{\theta_0, f}(x_i \mid X_n^i))^2}{V_{\theta_0}(x_i \mid X_n^i)} + \log \frac{V_{\theta_0}(x_i \mid X_n^i)}{ V_{\theta}(x_i \mid X_n^i)} \bigg] < 0,
    \end{equation} 
     then there exists $n_0 \in \R{N}$ such that $\cl(\theta_0 \mid X_n) < \cl(\theta \mid X_n)$ for every $\theta \in (0, \theta_1)$ and  $f \in \mathcal{B}$ for all $n \geq n_0$. Hence, $\inf_{f \in \mathcal{B}}\cle{\theta}(X_n) \geq \theta_1 $ when $n \geq n_0$.
    \end{enumerate}
\end{theorem}
\begin{proof}
    Given the condition~\eqref{gen_cond}, there exists $n_0 \in \R{N}$ such that 
    \[
    \sup_{\theta<\theta_1}\sup_{f \in \mathcal{B}}\bigg[ \| f \|^{2}_{\theta_0} + \log \frac{\det K_{\theta_0}(X_n)}{\det K_{\theta}(X_n)} \bigg] < 0
    \] 
    for all $n\geq n_0$. By Lemma~\ref{lem:LikelihoodComp}, this implies that $ \ml(\theta_0 \mid X_n) < \ml(\theta \mid X_n) $ for each $\theta \in (0, \theta_1)$ and for each $f \in \mathcal{B}$ when $n\geq n_0$. 
    Since $\mle{\theta}(X_n)$ is the argument that minimises $\ml(\theta \mid X_n)$, it follows that $\ml(\mle{\theta}(X_n) \mid X_n) < \ml(\theta \mid X_n)$ for every $\theta \in (0, \theta_1)$.
    Therefore, $\mle{\theta}(X_n)$ does not belong to $(0, \theta_1)$. That is, $\mle{\theta}(X_n) \geq \theta_1$. The proof for cross-validation is similar.    
\end{proof}

Hereafter, we consider the Gaussian kernel and its lengthscale.
We obtain lower bounds on the ML and CV estimators in two steps.
First, we use Lemma~\ref{lem:LikelihoodComp} to establish a condition on lengthscales $\alpha$ and $\beta$ that guarantees that $\ml(\alpha \mid X_n) < \ml(\beta \mid X_n)$ or $\cl(\alpha \mid X_n) < \cl(\beta \mid X_n)$ (Lemmas~\ref{th:ML_Compare} and~\ref{th:CV_Compare}).
The condition depends on the RKHS norm $\lVert f \rVert_\alpha$.
We then choose suitable classes of functions $\mathcal{B}$ that allow us to control $\lVert f \rVert_\alpha$ (Sections~\ref{sec:bounded} to~\ref{sec:power-exp}).
For the remainder of the article, we assume that the domain of interest is $\Omega = [-1, 1]$ and that the latent response function is the restriction onto $[-1,1]$ of $f \colon \mathbb{R} \to \mathbb{R}$.

\begin{lemma}\label{th:ML_Compare}
    Let $\delta>0$ and  $\alpha, \beta \geq \delta$. Suppose that $f \in \mathcal{H}_\alpha$.
    Then there exists a positive sequence $(\lambda_n)_n$ that converges to one such that $\ml(\alpha \mid X_n) < \ml(\beta \mid X_n)$ if 
     \[
     \beta <  \alpha \lambda_n \exp\bigg(-\frac{2\|f\|^{2}_{\alpha}}{n^2} \bigg)
     \]
    for any $n \geq 2$.  
\end{lemma}
\begin{proof}
By Lemma~\ref{lem:LikelihoodComp}, $\ml(\alpha \mid X_n) < \ml(\beta \mid X_n)$ if 
\begin{equation}\label{cond:ML}
     \| f \|^{2}_{\alpha} + \log \frac{\det K_{\alpha}(X_n)}{\det K_{\beta}(X_n)} < 0 .
\end{equation}
    The bounds on the Gaussian Gram determinant in Lemma~\ref{lem:DetBdFloater} bound the logarithmic term in~\eqref{cond:ML} as 
\begin{equation}\label{eq:LogRatioUpperBound}
    \log \frac{\det K_{\alpha}(X_n)}{\det K_{\beta}(X_n)} \leq n(n-1)\log\bigg(\frac{\beta}{\alpha}\bigg) + \frac{S_n}{\beta^2}.
\end{equation}
From the definition of $\Bar{x}$ and $S_n$ provided in~\eqref{eq:AverageMeanSquareDiff}, the condition $x_i, \Bar{x} \in [-1, 1]$ implies $|x_i - \Bar{x}| \leq 2$. This directly yields the upper bound $S_n \leq 4n$. Therefore 
\begin{equation*}\label{eq:Sandwich1}
    \| f \|^{2}_{\alpha} + \log \frac{\det K_{\alpha}(X_n)}{\det K_{\beta}(X_n)} \leq \| f \|^{2}_{\alpha} + n(n-1)\log\bigg(\frac{\beta}{\alpha}\bigg) + \frac{S_n}{\beta^2} \leq \| f \|^{2}_{\alpha} + n(n-1)\log\bigg(\frac{\beta}{\alpha}\bigg) + \frac{4n}{\delta^2}.
\end{equation*}
Dividing the right-hand side of this expression by $n(n-1)$, we obtain 
\begin{equation*}
    \frac{1}{n(n-1)}\| f \|^{2}_{\alpha} + \log\bigg(\frac{\beta}{\alpha}\bigg) + \frac{4}{(n-1)\delta^2}.
\end{equation*}
Using the inequality $n-1 \geq n/2$ (for $n \geq 2$), we can bound the terms involving $n-1$ to obtain
\begin{equation*}
     \frac{1}{n(n-1)}\| f \|^{2}_{\alpha} + \log\bigg(\frac{\beta}{\alpha}\bigg) + \frac{4}{(n-1)\delta^2} \leq  \frac{2}{n^2} \|f\|^{2}_{\alpha} + \log\bigg(\frac{\beta}{\alpha}\bigg) + \frac{8}{\delta^2 n}.
\end{equation*}
Consequently, if the condition
\begin{equation}\label{cond:AlmostFinalML_Comp}
     \frac{2}{n^2} \|f\|^{2}_{\alpha} + \log\bigg(\frac{\beta}{\alpha}\bigg) + \frac{8}{\delta^2 n} < 0,
\end{equation}
holds, then the condition in~\eqref{cond:ML} is satisfied. By exponentiating both sides, \eqref{cond:AlmostFinalML_Comp} can be equivalently expressed as 
\begin{equation}\label{cond:ML_Comp}
    \beta < \alpha \exp\bigg(-\frac{2\| f \|^{2}_{\alpha}}{n^2} - \frac{8}{\delta^2n}\bigg) = \alpha \lambda_n  \exp\bigg(-\frac{2\| f \|^{2}_{\alpha}}{n^2}\bigg), 
\end{equation}
where $ \lambda_n = e^{-8/(\delta^2n)}$ converges to $1$ as $n \to \infty$. 
It follows immediately that the relation $\ml(\alpha \mid X_n) < \ml(\beta \mid X_n)$ holds if $\beta$ and $\alpha$ satisfy the inequality~\eqref{cond:ML_Comp}.
\end{proof}
\begin{remark}
    If $\Omega = [a, b]$, then $S_n \leq n(b-a)^2 / 4$ by a form of Popoviciu's inequality~\cite{9a986506-935c-30f3-9482-ebde9bfc2680}.
    It follows that the above result can be generalised directly to any bounded subset of $\R{R}$.
\end{remark}

To establish a similar result for cross-validation we use Theorem~\ref{th:WCE_Bound}, which involves the constants $c_1(\ell)$ and $c_2(\ell)$.
In the following lemma, we derive a lower bound on the ratio of these constants.

\begin{lemma}\label{lower_bound_of_const}
Let $\delta>0$. Then there is a constant $C(\delta) > 0$ such that
    \[
      \frac{c_{1}(\beta)}{c_{2}(\alpha)} \geq \frac{C(\delta)}{\alpha}
     \] 
     for any $\alpha, \beta \geq \delta$.
\end{lemma}
\begin{proof}
    Since $x \mapsto \frac{x^2}{x^2+1}$ is strictly increasing on $(0, \infty)$ and $\beta \geq \delta > 0$, we have $\frac{\beta^2}{\beta^2+1} \geq \frac{\delta^2}{\delta^2+1}$. 
    Likewise, the exponent is increasing in $\beta$, so $-\frac{1}{32\beta^4} - \frac{1}{2\beta^2} \geq -\frac{1}{32\delta^4} - \frac{1}{2\delta^2}$. Thus, we can lower bound $c_1(\beta)$ by a constant depending only on $\delta$:
    \begin{equation*}
        c_{1}(\beta) = \sqrt{\frac{2\beta^2}{\beta^2 + 1}} \exp\bigg(-\frac{1}{32\beta^4} - \frac{1}{2\beta^2}\bigg)
        \geq \sqrt{\frac{2\delta^2}{\delta^2 + 1}} \exp\bigg(-\frac{1}{32\delta^4} - \frac{1}{2\delta^2}\bigg).
    \end{equation*}
We can upper bound $c_{2}(\alpha)$ by distributing the $\alpha^{1/2}$ term and factoring out $\alpha^2$. Since $\alpha \geq \delta$, we replace the remaining negative powers of $\alpha$ with negative powers of $\delta$:
\begin{align*}
    c_{2}(\alpha) &= \Bigg[\frac{\exp\big(\frac{1}{2\alpha^2}+ \frac{2}{\alpha}\big)}{2\sqrt{\pi} \alpha^{-1/2}}\bigg(3+\sqrt{\frac{\pi}{e}}+ 2\sqrt{\frac{\alpha}{\pi}}+ \frac{7}{16}\alpha + \frac{1}{\alpha}\Big(2+ \frac{1}{\alpha^2}\Big)\bigg)\Bigg]^{1/2}\\
    &= \Bigg[\frac{\exp\big(\frac{1}{2\alpha^2}+ \frac{2}{\alpha}\big)}{2\sqrt{\pi}}\alpha^2\bigg(3\alpha^{-3/2}+\sqrt{\frac{\pi}{e}}\alpha^{-3/2}+ \frac{2}{\sqrt{\pi}}\alpha^{-1}+ \frac{7}{16}\alpha^{-1/2} + 2\alpha^{-5/2} + \alpha^{-9/2}\bigg)\Bigg]^{1/2}\\
    &\leq \alpha \Bigg[\frac{\exp\big(\frac{1}{2\delta^2}+ \frac{2}{\delta}\big)}{2\sqrt{\pi}}\bigg(\frac{3}{\delta^{3/2}}+\sqrt{\frac{\pi}{e\delta^3}}+ 2\sqrt{\frac{1}{\pi\delta^2}}+ \frac{7}{16\sqrt{\delta}} + \frac{1}{\delta^{5/2}}\Big(2+ \frac{1}{\delta^2}\Big)\bigg)\Bigg]^{1/2}.
\end{align*}
Therefore, the lower bound of the ratio becomes
\begin{align*}
   \frac{c_{1}(\beta)}{c_{2}(\alpha)} \geq{}& \frac{1}{\alpha} \sqrt{\frac{2\delta^2}{\delta^2 + 1}} \exp\bigg(-\frac{1}{32\delta^4} - \frac{1}{2\delta^2} \bigg)\\
    &\times \bigg[ \frac{\exp\big(\frac{1}{2\delta^2}+ \frac{2}{\delta}\big)}{2\sqrt{\pi}}\bigg(\frac{3}{\delta^{3/2}}+\sqrt{\frac{\pi}{e\delta^3}}+ 2\sqrt{\frac{1}{\pi\delta^2}}+ \frac{7}{16\sqrt{\delta}} + \frac{1}{\delta^{5/2}}\Big(2+ \frac{1}{\delta^2}\Big)\bigg)\bigg]^{-1/2}\\
    \eqqcolon{} & \frac{C(\delta)}{\alpha}.  
\end{align*}
\end{proof}

\begin{lemma}\label{th:CV_Compare}
 Let $\delta>0$ and $\alpha, \beta \geq \delta$. 
 Suppose that $f \in \mathcal{H}_\alpha$.
 Then there exists a positive sequence $(\lambda_n)_{n}$ converging to one such that $\cl(\alpha \mid X_n) < \cl(\beta \mid X_n)$ if 
     \begin{equation}\label{cond:CV_Final}
     4\beta < \alpha^{1-\frac{1}{n-1}}\,\lambda_n \exp\bigg({-\frac{\|f\|^{2}_{\alpha}}{n}} \bigg)
     \end{equation}
 for any $n \geq 2$.
\end{lemma}
\begin{proof}
By Lemma~\ref{lem:LikelihoodComp}, $\cl(\alpha \mid X_n) < \cl(\beta \mid X_n)$ if 
\begin{equation}\label{cond:CV}
       \sum_{i=1}^{n} \bigg[ \frac{(f(x_i) - \mu_{\alpha, f}(x_i \mid X_n^i))^2}{V_{\alpha}(x_i \mid X_n^i)} + \log \frac{V_{\alpha}(x_i \mid X_n^i)}{ V_{\beta}(x_i \mid X_n^i)} \bigg] < 0 .
    \end{equation}
    We proceed as in the proof of Lemma~\ref{th:ML_Compare}.
    By applying the bounds on $V_{\ell}(x \mid X_n)$ from Theorem~\ref{th:WCE_Bound}, we obtain
    \[
    \frac{V_{\alpha}(x_i \mid X_n^i)}{ V_{\beta}(x_i \mid X_n^i)} \leq \bigg(
    \frac{c_{2}(\alpha)}{c_{1}(\beta)}\bigg)^2\bigg( \frac{4 \beta}{\alpha} \bigg)^{2(n-1)} (n-1)^{-1/4} e^{2\sqrt{n-1}/\alpha}.
    \]
    Taking the natural logarithm on the above inequality and summing yields
    \[
    \sum_{i=1}^n\log \frac{V_{\alpha}(x_i \mid X_n^i)}{ V_{\beta}(x_i \mid X_n^i)}  \leq 2n\log \bigg(
    \frac{c_{2}(\alpha)}{c_{1}(\beta)}\bigg) -\frac{n}{4} \log(n-1) + 2n(n-1)\log \bigg( \frac{4 \beta}{\alpha} \bigg) + \frac{2n^{3/2}}{\alpha}.
    \]
    From Lemma~\ref{lower_bound_of_const} we know that $c_1(\beta)/ c_2(\alpha) \geq C(\delta)/\alpha$, which implies $c_2(\alpha)/c_1(\beta) \leq \alpha/C(\delta)$.
    From this, the estimate $\log(n-1) \geq \log(n/2)$ (which holds when $n \geq 2$), and the fact that $\alpha \geq \delta$, we obtain    
    \[
    \sum_{i=1}^n\log \bigg(\frac{V_{\alpha}(x_i \mid X_n^i)}{ V_{\beta}(x_i \mid X_n^i)}\bigg) \leq 2n\log \bigg( \frac{\alpha}{C(\delta)}\bigg) -\frac{n}{4} \log \bigg(\frac{n}{2}\bigg) + 2n(n-1)\log \bigg( \frac{4 \beta}{\alpha} \bigg) + \frac{2n^{3/2}}{\delta}.
    \]
    With \eqref{eq:CVNorm} the first term in~\eqref{cond:CV} is bounded as    
    \[\sum_{i=1}^n \frac{(f(x_i) - \mu_{\alpha, f}(x_i \mid X_n^i))^2}{V_{\alpha}(x_i \mid X_n^i)} \leq n \|f\|_{\alpha}^2.\]    
    Consequently, the left-hand side of~\eqref{cond:CV} is bounded from above by
    \[
    n \|f\|_{\alpha}^2 + 2n\log \bigg( \frac{\alpha}{C(\delta)}\bigg) -\frac{n}{4} \log \bigg(\frac{n}{2}\bigg) + 2n(n-1)\log \bigg( \frac{4 \beta}{\alpha} \bigg) + \frac{2n^{3/2}}{\delta}.
    \] 
    Dividing this expression by $2n(n-1)$, we obtain
    \[
    \frac{1}{2(n-1)} \|f\|_{\alpha}^2 + \frac{1}{n-1}\log \bigg( \frac{\alpha}{C(\delta)}\bigg) -\frac{1}{8(n-1)} \log \bigg(\frac{n}{2}\bigg) + \log \bigg( \frac{4 \beta}{\alpha} \bigg) + \frac{n^{1/2}}{(n-1)\delta}.
    \]
    Using again the inequality $n-1 \geq n/2$, we can bound the terms involving $n-1$ to obtain
    \begin{align}\label{eq:IntermedIneq3.6}
    \frac{1}{2(n-1)} \|f\|_{\alpha}^2 &+ \frac{1}{n-1}\log \bigg( \frac{\alpha}{C(\delta)}\bigg) -\frac{1}{8(n-1)} \log \bigg(\frac{n}{2}\bigg) + \log \bigg( \frac{4 \beta}{\alpha} \bigg) + \frac{n^{1/2}}{(n-1)\delta} \notag \\
        &\leq \frac{\|f\|_{\alpha}^2}{n} + \frac{1}{n-1}\log \bigg( \frac{\alpha}{C(\delta)}\bigg) -\frac{1}{8n} \log \bigg(\frac{n}{2}\bigg) + \log \bigg( \frac{4 \beta}{\alpha} \bigg) + \frac{2n^{-1/2}}{\delta} < 0, 
    \end{align}
    which implies~\eqref{cond:CV}.
    Rearranging the last inequality in~\eqref{eq:IntermedIneq3.6} gives
    \begin{equation}\label{cond:CV-Comp}
     4 \beta <  \alpha^{1-\frac{1}{n-1}} C(\delta)^{\frac{1}{n-1}} \bigg( \frac{n}{2} \bigg)^{\frac{1}{8n}} e^{-\frac{2}{\delta\sqrt{n}}} \exp \bigg( - \frac{\|f\|_{\alpha}^2}{n} \bigg) = \alpha^{1-\frac{1}{n-1}} \lambda_n \exp \bigg( - \frac{\|f\|_{\alpha}^2}{n} \bigg),
    \end{equation}
    where
    \[
    \lambda_n = C(\delta)^{\frac{1}{n-1}} \bigg( \frac{n}{2} \bigg)^{\frac{1}{8n}} e^{-\frac{2}{\delta\sqrt{n}}}.
    \]     
    Therefore, $\cl(\alpha \mid X_n)<\cl(\beta \mid X_n)$ if~\eqref{cond:CV-Comp} holds. 
\end{proof}
If $f \in \mathcal{H}_\ell$, Theorem~\ref{th:BdSubsetRKHS} below shows that $\liminf_{n \to \infty} \mle{\ell}(X_n) \geq \ell$ but $\liminf_{n \to \infty} \cle{\ell}(X_n) \geq \ell/4$.
The difference stems from the factor $4$ in~\eqref{cond:CV_Final}, which is inherited from the lower bound on $V_\ell(x\mid X_n)$ given in Theorem~\ref{th:WCE_Bound}.     
One obvious option to remove this factor is to use Gaussian Gram determinant bound provided in Lemma~\ref{lem:DetBdFloater} and the identity~\eqref{eq:VarToDet}.
    In this way one can show that $\cl(\alpha \mid X_n) < \cl(\beta \mid X_n)$ if 
     \[
     \beta < \alpha \,e^{-{q_2}/{\delta^2}}\exp\bigg({-\frac{\|f\|^{2}_{\alpha}}{n}} \bigg) 
     \]
     for some positive constant $q_2$.
     While selecting a large $\delta$ makes the factor $e^{-{q_2}/{\delta^2}}$ close to one, which is desirable, it also restricts the parameter space to $\Theta = [\delta, \infty)$, which does not necessarily include the ``correct'' value of the parameter.
     

Next we impose various assumptions on the latent function $f$.
These assumptions allow us to control the norm $\lVert f \rVert_\alpha$ in~\eqref{cond:ML} and~\eqref{cond:CV_Final}.

\subsection{Bounded subset of $\mathcal{H}_{\ell}$} \label{sec:bounded}

 Let $\delta >0$ and define our parameter set for lengthscale as $\Theta = [\delta, \infty)$. 
 That is, we restrict the estimators to lie in the interval $[\delta, \infty)$.
 Here we assume that the latent function is an element of the RKHS $\mathcal{H}_\ell$ for some $\ell >0$. The subsequent theorem provides a lower bound for CV and ML estimators of the lengthscale.

\begin{theorem}\label{th:BdSubsetRKHS}
    Let $\Theta = [\delta, \infty)$ for some $\delta > 0$.
    If $\mathcal{B}$ is a bounded subset of $\mathcal{H}_{{\ell}}$, then \[
     \liminf_{n \to \infty}\inf_{f \in \mathcal{B}}\mle{\ell}(X_n) \geq \ell
     \quad \text{and} \quad \liminf_{n \to \infty}\inf_{f \in \mathcal{B}}\cle{\ell}(X_n) \geq \frac{\ell}{4}.
     \]
\end{theorem}
\begin{proof}
     The case $\ell < \delta$ is trivial because both estimators are greater than or equal to $\delta$. We consider then the case $\ell \geq \delta$. By Lemma~\ref{th:ML_Compare}, for $n \geq 2$ and $f \in \mathcal{B}$ we have $\ml(\ell \mid X_n) < \ml(\beta \mid X_n)$ if
     \[
     \beta < \ell \lambda_n
     \exp{\bigg(-\frac{2\|f\|^{2}_{\ell}}{n^2}\bigg)},
     \] 
     where $\beta \geq \delta$ and $(\lambda_n)_n$ converges to $1$. 
     By Definition~\ref{mle_def}, this implies 
     \[
     \mle{\ell}(X_n) \geq \ell \lambda_n
     \exp{(-{2\|f\|^{2}_{\ell}}/{n^2})}
     \] 
     for any $f \in \mathcal{B}$. 
     Since $\mathcal{B}$ is a bounded subset of $\mathcal{H}_{{\ell}}$, there exists a constant $m>0$ such that $\|f\|_\ell \leq m$ for all $f \in \mathcal{B}$. Consequently,
    \[
    \inf_{f \in \mathcal{B}}\exp{\bigg(-\frac{2\|f\|^{2}_{\ell}}{n^2}\bigg)} \geq \exp{\bigg(-\frac{2m^2}{n^2}\bigg)}.
    \]
     For each $f \in \mathcal{B}$, define the sequence $(a_n^f)_n$ by 
    \[
    a_n^f = {\ell} \lambda_n \exp{\bigg(-\frac{2\|f\|^{2}_{\ell}}{n^2}\bigg)}.
    \]
    Then $\ml(\ell \mid X_n ) < \ml(\beta \mid X_n)$ for each $\beta \in [\delta, a_n^f)$. This implies $\inf_{f \in \mathcal{B}}\mle{\ell}(X_n) \geq \inf_{f \in \mathcal{B}} a_n^f$. Since $(\lambda_n)_n$ depends only on $\delta$ and $\lim_{n \to \infty} \lambda_n = 1$, we have   
    \[
    \liminf_{n\to \infty} \inf_{f \in \mathcal{B}}\mle{\ell}(X_n) \geq \liminf_{n \to \infty} \inf_{f \in \mathcal{B}} a_n^f \geq \lim_{n \to \infty} {\ell}\lambda_n \exp{\bigg(-\frac{2m^2}{n^2}\bigg)} = {\ell}. 
    \] 
    The proof for lower bound of CV estimator is analogous and involves replacing Lemma~\ref{th:ML_Compare} with Lemma~\ref{th:CV_Compare}. 
\end{proof}

Let $f \in \mathcal{H}_{\ell}$ for all $\ell \in (0, \ell_0)$ but $f \notin \mathcal{H}_{\ell}$ for each $\ell > \ell_0$ and assume that $\delta < \ell_0$.
Then Theorem~\ref{th:BdSubsetRKHS} yields 
\begin{equation*}
    \liminf_{n \to \infty} \mle{\ell}(X_n)  \geq \sup_{\ell \in [\delta, \ell_0)} \ell =\ell_0 \quad \text{ and } \quad \liminf_{n \to \infty} \cle{\ell}(X_n)  \geq \ell_0/4.
\end{equation*}
See Theorem~\ref{th:Minh2010} for a result regarding function inclusion in the RKHS of the Gaussian kernel.

\subsection{Bandlimited functions}

Next we assume that the latent response function is bandlimited.
This means that its Fourier transform is compactly supported.
For such functions, the ML and CV estimators diverge at least as fast as $\sqrt{\log n}$.

\begin{theorem}\label{th:BandLimit}
    Let $\Theta = [\delta, \infty)$ for some $\delta > 0$.
    Let $0<c_0<1$ and $m, R>0$. Define
    \[
    \mathcal{B}_R^m = \big\{ f \colon \R{R}\to \R{R} \: | \: \supp(\hat{f})\subset [-R,R]\quad \text{and} \quad  \|\hat{f}\|_{L^2(\R{R})} < m \big\} .
    \]
    Then 
    \begin{equation*}
        \inf_{f \in \mathcal{B}_R^m}\mle{\ell}(X_n) \geq \bigg(\frac{c_0}{R}\bigg) \sqrt{\log n+ \kappa_1 } \quad \text{and} \quad \inf_{f \in \mathcal{B}_R^m} \cle{\ell}(X_n) \geq \bigg(\frac{\sqrt{2}\,c_0}{16\,R}\bigg)
    \sqrt{\log n + \kappa_2 }
    \end{equation*}
    for sufficiently large $n$ with 
    constants $\kappa_1$ and $\kappa_2$ defined in~\eqref{eq:kappa1} and~\eqref{eq:kappa2}, respectively.
    
\end{theorem}
\begin{proof}
    Consider the ML estimator.
    By Lemma~\ref{th:ML_Compare}, $\ml(\alpha \mid X_n) < \ml(\beta \mid X_n)$ if 
     \[
     \beta < \alpha\lambda_n \exp\bigg(-\frac{2\|f\|^{2}_{\alpha}}{n^2} \bigg),
     \]
     where $(\lambda_n)_{n}$ converges to 1 and $\alpha, \beta \geq \delta$.
     Let $A_n(\alpha) =  \alpha\lambda_n \exp\big({-{2\|f\|^{2}_{\alpha}}/{n^2}}\big)$. 
     Then, by the definition of MLE we infer $\mle{\ell}(X_n) \geq A_n(\alpha)$.  
     Since we are seeking the sharpest possible lower bound for $\mle{\ell}(X_n)$ using the theory developed herein, our goal is to find a sequence $(\alpha_n)_n$ that maximises the rate of growth of $A_n(\alpha_n)$ as $n$ tends to infinity. 
     Let us compute an upper bound of the norm of $f \in \mathcal{B}_R^m$:
    \begin{align*}
        \|f\|^2_{\alpha} = \frac{1}{\sqrt{2\pi}\alpha}\int_{-R}^{R}|\hat{f}(\omega)|^2 e^{\alpha^2\omega^2/2}\,d\omega
        &\leq \frac{e^{\alpha^2R^2/2}}{\sqrt{2\pi}\alpha}\|\hat{f}\|_{L^2(\R{R})}^2\\
        &\leq \frac{e^{\alpha^2R^2/2}}{\sqrt{2\pi}\delta}m^2.
    \end{align*}
    Since $\|f\|_\alpha^2$ grows at a rate not greater than $e^{(\alpha^2R^2/2)}$ in $\alpha$, in order to achieve the fastest possible growth rate of $A_n(\alpha)$, we require $\exp({-{2\|f\|^{2}_{\alpha}}/{n^2}})\geq c_0$ for some constant $0< c_0 <1$, or equivalently $2\|f\|_{\alpha}^2 \leq -n^2 
     \log c_0$.  
    Therefore, $2\|f\|_{\alpha}^2 \leq -n^2\log c_0$ holds if 
    \begin{align}\label{eq:norm_bound}
        \|f\|_{\alpha}^2 \leq \frac{e^{\alpha^2R^2/2}m^2}{\sqrt{2\pi}\delta} \leq \frac{n^2}{2}
        \log\bigg(\frac{1}{c_0}\bigg).
    \end{align}
     By considering the second inequality from \eqref{eq:norm_bound}, we obtain an upper bound on $\alpha$ such that the above inequality holds:
    \begin{equation}\label{eq:UpperBdAlphaML}
        \alpha^2  \leq \frac{2}{R^2} \log\bigg(n^2\frac{\sqrt{2\pi}\delta\log (1/c_0) }{2m^2}\bigg).
    \end{equation}
    Thus, for each $n$, we choose $\alpha_n$ as the largest admissible $\alpha$ from~\eqref{eq:UpperBdAlphaML} to ensure that $\exp({-{2\|f\|^{2}_{\alpha}}/{n^2}}) \geq c_0$. This yields
    \[
    \alpha_n= \frac{2}{R}(\log n + \kappa_1 )^\frac{1}{2},
    \]
     where 
    \begin{equation}\label{eq:kappa1}
        \kappa_1 = \frac{1}{2}\log\bigg[\frac{\sqrt{2\pi} \delta\log (1/c_0) }{2m^2}\bigg].
    \end{equation}
    Then Lemma~\ref{th:ML_Compare} states that $\ml(\alpha_n \mid X_n) < \ml(\beta \mid  X_n)$ if 
    \[
    \beta < c_0\lambda_n\,\alpha_n 
    \]
    for sufficiently large $n$, where $\beta \geq \delta$. Since $(c_0\lambda_n)_n$ converges to $c_0$, $c_0/2<c_0 \lambda_n$ for any sufficiently large $n$.
    Therefore $\mle{\ell}(X_n)>  c_0\,\alpha_n/2$ for any sufficiently large $n$ by Definition~\ref{mle_def}. 

    For the proof of the lower bound of $\cle{\ell}(X_n)$, we follow the same technique as for $\mle{\ell}$. The proof involves replacing Lemma~\ref{th:ML_Compare} with Lemma~\ref{th:CV_Compare} and introducing the condition $\|f\|_{\alpha}^2 \leq -n\log c_0$. 
    This allows us to select a sequence $(\beta_n)_n$ similar to $(\alpha_n)_n$ such that
    \[
    \beta_n = \frac{\sqrt{2}}{R} \big[\log n + \kappa_2 \big]^{\frac{1}{2}},
    \] 
    where 
    \begin{equation}\label{eq:kappa2}
        \kappa_2 = \log\bigg[ \frac{\sqrt{2\pi} \delta \log (1/c_0)}{m^2}\bigg] = 2\kappa_1+\log2. 
    \end{equation}
Therefore, for any $ f \in \mathcal{B}_R^m$,
\[
\cle{\ell}(X_n) \geq \Big(\frac{c_0}{8}\Big)\, \beta_n^{(1-\frac{1}{n-1})} = \Big(\frac{c_0}{8}\Big)\, \bigg[\frac{\sqrt{2}}{R} \sqrt{\log n + \kappa_2 }\bigg]^{(1-\frac{1}{n-1})}.
\]
Since $\smash{(({\sqrt{2}}/{R}) \sqrt{\log n + \kappa_2 })^{{1}/({n-1})}}$ converges to one as $n$ tends to $\infty$, for all sufficiently large $n$, 
\[
\bigg[\bigg(\frac{\sqrt{2}}{R}\bigg) \sqrt{\log n + \kappa_2 }\bigg]^{{1}/({n-1})} \leq 2.
\]
Hence, for any $f \in \mathcal{B}_R^m$,
    \[
    \cle{\ell}(X_n) \geq \bigg(\frac{\sqrt{2}\,c_0}{16\,R}\bigg)
    \sqrt{\log n + \kappa_2 }.
    \]    
\end{proof}

\subsection{Power exponential Fourier transform} \label{sec:power-exp}

 Here we consider functions $f$ whose Fourier transforms satisfy $\lvert \hat{f}(\omega) \rvert \leq M \exp(-a|\omega|^\gamma)$ for $\gamma >2$ and $a, M >0$. 
 Such functions lie in $\mathcal{H}_\ell$ because
\begin{align*}
    \|f\|^2_\ell = \frac{1}{\sqrt{2\pi}\ell} \int_{\R{R}} |\hat{f}(\omega)|^2 \exp\bigg( \frac{\ell^2 \omega^2}{2} \bigg) d\omega 
     & \leq \frac{1}{\sqrt{2\pi}\ell} \int_{\R{R}}M^2 \exp(-2a|\omega|^\gamma) \exp\bigg( \frac{\ell^2 \omega^2}{2} \bigg) d\omega\\ 
     &= \frac{M^2}{\sqrt{2\pi}\ell} \int_{\R{R}} \exp\bigg(-2a|\omega|^\gamma + \frac{\ell^2 \omega^2}{2} \bigg) d\omega,
\end{align*}
which is obviously finite if $\gamma>2$.

\begin{theorem}\label{th:PowerExpo}
    Let $\Theta = [\delta, \infty)$ for some $\delta > 0$.
    Let $\gamma, a, M > 0$ and $0<c_0<1$.
    Define \[
    \mathcal{B}_\gamma = \big\{f\colon \R{R} \to \R{R} \colon |\hat{f}(\omega)| \leq M \exp(-a|\omega|^\gamma) \big\} .
    \]     
    If $\gamma > 2$, then 
    \begin{align*}
    \inf_{f \in \mathcal{B}_\gamma}\mle{\ell}(X_n) 
    &\geq \bigg(\frac{c_0}{2}\bigg) \bigg(\frac{2}{R}\bigg)^{\frac{\gamma-2}{\gamma}} (\log n + \rho_1)^{\frac{\gamma-2}{2\gamma}} 
    \end{align*}
    and
    \begin{align*}
     \inf_{f \in \mathcal{B}_\gamma}\cle{\ell}(X_n) \geq \bigg(\frac{c_0}{16}\bigg) \bigg(\frac{\sqrt{2}}{R}\bigg)^{\frac{\gamma-2}{\gamma}} (\log n + \rho_2)^{\frac{\gamma-2}{2\gamma}} 
    \end{align*}
        for sufficiently large $n$ with $R = (2a)^{-1/(\gamma-2)}$ and constants $\rho_1$ and $\rho_2$  defined in \eqref{eq:rho_1} and \eqref{eq:rho_2}, respectively.
\end{theorem}
\begin{proof}
    The proof is analogous to that of Theorem~\ref{th:BandLimit}. Specifically, we utilise the condition $2\|f\|_{\alpha}^2 \leq -n^2 \log c_0$ to select a sequence $(\alpha_n)_n$ (or $(\beta_n)_n$ in the case of CV estimator) that establishes the proposed lower bound for the ML estimator. To achieve this, we employ Lemma~\ref{th:ML_Compare} (alternatively Lemma~\ref{th:CV_Compare} for CV estimator) alongside the properties of Fourier transform of functions in $\mathcal{B}_\gamma$. 
    
    We have
    \begin{align}\label{int:ProImpro}
        \|f\|^2_\alpha &= \frac{1}{\sqrt{2\pi}\alpha} \int_{-\infty}^{\infty} |\hat{f}(\omega)|^2 e^{\alpha^2\omega^2/2}\,d\omega \notag\\
        & \leq M^2\bigg[ \frac{2}{\sqrt{2\pi}\alpha} \int_{0}^{r} e^{-2a\omega^\gamma} e^{\alpha^2\omega^2/2}\,d\omega +  \frac{2}{\sqrt{2\pi}\alpha} \int_{r}^{\infty}e^{-2a\omega^\gamma} e^{\alpha^2\omega^2/2}\,d\omega \bigg]
    \end{align}
    for any $r > 0$.
    Let us estimate the two integrals. 
   Set $r = (\alpha^2/(2a))^{1/(\gamma-2)}$.
   Then $-2a\omega^\gamma + \alpha^2 \omega^2 / 2 \leq -a\omega^\gamma$ if $\omega \geq r$.
    Therefore
   \begin{equation}\label{int:Improper}  \int_{r}^{\infty} e^{-2a\omega^\gamma} e^{\alpha^2\omega^2/2} \,d\omega \leq \int_{0}^{\infty} e^{-a\omega^\gamma}\,d\omega = \frac{1}{\gamma}a^{-1/\gamma} \Gamma \bigg(\frac{1}{\gamma}\bigg), 
   \end{equation}
   where the substitution $x = a^{1/\gamma} \omega$ and the definition of the gamma function were used to compute the integral.
    The proper integral in~\eqref{int:ProImpro} can be estimated as
    \begin{equation}\label{int:Proper}
        \int_{0}^{r} e^{-2a\omega^\gamma} e^{\alpha^2\omega^2/2}\,d\omega \leq \int_{0}^{r} e^{\alpha^2\omega^2/2}\,d\omega
        \leq \int_{0}^{r} e^{\alpha^2r\omega/2}\,d\omega 
        = \frac{2}{\alpha^2r}\big(e^{\alpha^2r^2/2}-1\big) .
    \end{equation}
    Combining the upper bounds in~\eqref{int:Improper} and~\eqref{int:Proper} yields
    \begin{align*}
        \|f\|^2_\alpha &\leq \frac{2M^2}{\sqrt{2\pi}\alpha}\bigg( \frac{1}{\gamma}a^{-1/\gamma}\Gamma\bigg(\frac{1}{\gamma}\bigg) + \frac{2}{\alpha^2r}\big(e^{\alpha^2r^2/2}-1\big)\bigg) \\ 
        & \leq \frac{2M^2}{\sqrt{2\pi}\alpha}\bigg( \frac{1}{\gamma}a^{-1/\gamma}\Gamma\bigg(\frac{1}{\gamma}\bigg) + \frac{2}{R\delta^{2(\gamma - 1) / (\gamma-2)}}\big(e^{\alpha^2r^2/2}-1\big)\bigg).
    \end{align*} 
    By substituting the previously defined value of $r$ into the expression, and defining the constants $R = (2a)^{-1/(\gamma-2)}$, $d = \delta^{2(\gamma - 1) / (\gamma-2)}$, and $b = \frac{1}{\gamma}a^{-1/\gamma}\Gamma({1}/{\gamma})$, we obtain
    \[
    \|f\|^2_\alpha 
    \leq \frac{2M^2}{\sqrt{2\pi}\delta} \bigg(  b + \frac{2}{Rd}\exp\bigg(\frac{R^2}{2}\alpha^{2(1+\frac{2}{\gamma-2})}\bigg)\bigg) .
    \]
    Since $2\|f\|^2_\alpha$ has a growth rate analogous to that in the proof of Theorem~\ref{th:BandLimit}, we can apply a similar argument to obtain an upper bound for $\alpha$.
    That is, we use the condition 
    \begin{equation}\label{eq:NormLessThann^2logC}
    2\|f\|^2_\alpha \leq \frac{4M^2}{\sqrt{2\pi}\delta} \bigg(  b + \frac{2}{Rd}\exp\bigg(\frac{R^2}{2}\alpha^{\frac{2\gamma}{\gamma-2}}\bigg)\bigg)\leq -n^2\log c_0.
    \end{equation}
    Rearranging this inequality gives
    \begin{align*}
    \alpha^{\frac{2\gamma}{\gamma-2}} &\leq \frac{2}{R^2}\bigg(\log\bigg(\sqrt{\frac{\pi}{2}} \frac{n^2\delta}{2M^2} \log \bigg(\frac{1}{c_0}\bigg) -b\bigg) +\log\bigg(\frac{Rd}{2}\bigg) \bigg) \notag \\
    &= \bigg(\frac{2}{R}\bigg)^2\bigg(\log n + \frac{1}{2}\log\bigg(\sqrt{\frac{\pi}{2}} \frac{\delta}{2M^2} \log \bigg(\frac{1}{c_0}\bigg) -\frac{b}{n^2}\bigg) +\frac{1}{2}\log\bigg(\frac{Rd}{2}\bigg) \bigg).
    \end{align*}
    Since $\sqrt{\frac{\pi}{2}} \frac{\delta}{2M^2} \log \big(\frac{1}{c_0}\big) -\frac{b}{n^2}$ converges to $\sqrt{\frac{\pi}{2}} \frac{\delta}{2M^2} \log \big(\frac{1}{c_0}\big)$ as $n \to \infty$, 
    \[
    \sqrt{\frac{\pi}{2}} \frac{\delta}{2M^2} \log \bigg(\frac{1}{c_0}\bigg) -\frac{b}{n^2} \geq \frac{\sqrt{\pi}}{4\sqrt{2}} \frac{\delta}{M^2} \log \bigg(\frac{1}{c_0}\bigg)
    \]
    for all sufficiently large $n$. This implies
    \begin{align*}
    \bigg(\frac{2}{R}\bigg)^2& \bigg(\log n + \frac{1}{2}\log\bigg(\frac{\sqrt{\pi}}{4\sqrt{2}} \frac{\delta}{M^2} \log \bigg(\frac{1}{c_0}\bigg)\bigg) +\frac{1}{2}\log\bigg(\frac{Rd}{2}\bigg) \bigg) \\
    &\leq \bigg(\frac{2}{R}\bigg)^2 \bigg(\log n + \frac{1}{2}\log\bigg(\sqrt{\frac{\pi}{2}} \frac{\delta}{2M^2} \log \bigg(\frac{1}{c_0}\bigg)-\frac{b}{n^2}\bigg) +\frac{1}{2}\log\bigg(\frac{Rd}{2}\bigg) \bigg).
    \end{align*}
    For each $n$, we choose $\alpha_n$ such that
    \begin{align*}
    \alpha_n^{\frac{2\gamma}{\gamma-2}} 
    &= \bigg(\frac{2}{R}\bigg)^2\bigg(\log n + \frac{1}{2}\log\bigg(\frac{\sqrt{\pi}}{4\sqrt{2}} \frac{\delta}{M^2}\log \bigg(\frac{1}{c_0}\bigg)\bigg) +\frac{1}{2}\log\bigg(\frac{Rd}{2}\bigg) \bigg) \notag \\
    &= \bigg(\frac{2}{R}\bigg)^2\big[\log n + \rho_1 \big] 
    \end{align*}
    where \begin{equation}\label{eq:rho_1}
    \rho_1 \coloneqq \frac{1}{2}\log\bigg(\frac{\sqrt{2\pi} \delta \log({1}/{c_0})}{8M^2}\bigg) +\frac{1}{2}\log\bigg(\frac{Rd}{2}\bigg).
    \end{equation}
    Hence $2\|f\|_{\alpha_n}^2 \leq -n^2 \log c_0$ holds for all sufficiently large $n$.
    Therefore, applying Lemma~\ref{th:ML_Compare}, we obtain  $\ml(\alpha_n \mid X_n) < \ml(\beta \mid  X_n)$ if 
    \[
    \beta < c_0\lambda_n\;\alpha_n
    \]
    for sufficiently large $n$, where $\beta \geq \delta$. Since $(c_0\lambda_n)_n$ converges to $c_0$, $c_0/2<c_0 \lambda_n$ for any sufficiently large $n$. Therefore $\mle{\ell}(X_n)>  c_0\,\alpha_n/2$ for sufficiently large $n$ by Definition~\ref{mle_def}.

    For the proof of the lower bound of $\cle{\ell}(X_n)$, we follow the same technique as for $\mle{\ell}$. The proof involves replacing Lemma~\ref{th:ML_Compare} with Lemma~\ref{th:CV_Compare} and introducing the condition $\|f\|_{\alpha}^2 \leq -n\log c_0$. This allows us to select a sequence $(\beta_n)_n$ analogous to $(\alpha_n)_n$ such that
    \[
    \beta_n^{\frac{2\gamma}{\gamma -2}} = \frac{2}{R^2} (\log n + \rho_2 ),
    \] 
    where 
    \begin{equation}\label{eq:rho_2}
    \rho_2 = 2\rho_1 + \log 2.
    \end{equation} Therefore, for any $f \in \mathcal{B}_\gamma$, 
    \begin{align*}
    \cle{\ell}(X_n) \geq \Big(\frac{c_0}{8}\Big)\, \beta_n^{(1-\frac{1}{n-1})} &= \Big(\frac{c_0}{8}\Big)\, \bigg(\bigg[\frac{2}{R^2} (\log n + \rho_2 )\bigg]^{(1-\frac{1}{n-1})} \bigg)^{\frac{\gamma -2}{2\gamma}}\\
    & = \Big(\frac{c_0}{8}\Big)\bigg[\frac{\sqrt{2}}{R} \sqrt{\log n + \rho_2 }\bigg]^{\frac{\gamma -2}{\gamma}}\bigg(\bigg[\frac{\sqrt{2}}{R} \sqrt{\log n + \rho_2}\bigg]^{-\frac{1}{n-1}} \bigg)^{\frac{\gamma -2}{\gamma}}
    \end{align*}
    for sufficiently large $n$. Since $\smash{(({\sqrt{2}}/{R}) \sqrt{\log n + \rho_2 })^{{1}/({n-1})}}$ converges to one as $n$ tends to $\infty$, for all sufficiently large $n$,
    \[
    \bigg[\bigg(\frac{\sqrt{2}}{R}\bigg) \sqrt{\log n + \rho_2 }\bigg]^{{1}/({n-1})} \leq 2.
    \]
    Hence, for any $f \in \mathcal{B}_\gamma$, 
    \[
    \cle{\ell}(X_n) \geq \bigg(\frac{c_0}{8}\bigg)\,\bigg[\frac{\sqrt{2}}{2 R} \sqrt{\log n + \rho_2} \bigg]^{\frac{\gamma -2}{\gamma}} \geq \bigg(\frac{c_0}{16}\bigg)\,\bigg[\frac{\sqrt{2}}{R} \sqrt{\log n + \rho_2 }\bigg]^{\frac{\gamma -2}{\gamma}}.
    \]
\end{proof}
\begin{remark}    
    Informally, bandlimited functions with bounded Fourier transform can be interpreted as functions with $\gamma = \infty$: if $f$ is bandlimited with bounded Fourier transform, then for every $\gamma > 2$ there is $a, M > 0$ such that $f \in \mathcal{B}_\gamma$.
    The rates in Theorem~\ref{th:PowerExpo} for $\gamma = \infty$ coincide with those in Theorem~\ref{th:BandLimit}.
\end{remark}
 

\section{Numerical experiments}

In this section, we present numerical experiments demonstrating the asymptotic lower bound of the lengthscale provided in Theorem~\ref{th:BdSubsetRKHS}. 
Due to ill-conditioning of the Gaussian Gram matrix and the polylogarithmic rate of growth of the lengthscale estimators, it proved numerically challenging to illustrate Theorems~\ref{th:BandLimit} and~\ref{th:PowerExpo}.

Let $k_{x_0}^\ell(x) = \exp\big(-{(x-x_0)^2}/{2\ell^2}\big)$ for $x_0 \in \mathbb{R}$. 
In the experiments we assume that the response function is a linear combination of such functions.
The following theorem by Minh~\cite{minh2010some} characterises when $k_{x_0}^{\ell'}$ is contained in the RKHS $\mathcal{H}_\ell$.

\begin{theorem}[Theorem 3 in \cite{minh2010some}]\label{th:Minh2010}
Let $\Omega \subset \mathbb{R}$ be any set with non-empty interior and $x_0 \in \mathbb{R}$. 
Let $K_\ell(x,y) = \exp\big(-{(x-y)^2}/{2\ell^2}\big)$. 

Then the Hilbert space $\mathcal{H}(K_\ell, \Omega)$ induced by $K_\ell$ on $\Omega$ contains the function $k^{\ell/\sqrt{\mu}}_{x_0}$ 
if and only if $0 < \mu < 2$. 
\end{theorem}

Fix $\ell >0$. Let $\nu \in (0,2)$ and $\ell_0 = \sqrt{\nu}\ell$. Then   $K_{\ell_0}(x,y) = \exp\left(-{(x-y)^2}/{2\nu\ell^2}\right)$ and 
\[
 k_{x_0}^{\ell\sqrt{\nu/\mu}} \in \mathcal{H}(K_{\ell_0}, [-1, 1]) = \mathcal{H}_{\ell_0}
\] 
if and only if $0<\mu<2$. If we choose $\mu = \nu$, then $k_{x_0}^{\ell} \in \mathcal{H}_{\ell_0}$ for any $0<\ell_0<\sqrt{2}\ell$. Let $\ell_1<\ell_2<\dots<\ell_m$ and
\begin{equation}\label{eq:LinearSumofFeatures}
f(x)=\sum_{i=1}^m c_i k_{z_i}^{\ell_i}(x)
\end{equation}
for any $c_i \neq 0, i = 1,2,\dots, m$. Then $f \in \mathcal{H}_{\ell_0}$ for all $\ell_0 \in (0,\sqrt{2}\ell_1)$ and $f \notin \mathcal{H}_{\ell_0}$ for $\ell_0 \geq \sqrt{2} \ell_1$.
Therefore, by Theorem~\ref{th:BdSubsetRKHS}, we have $\liminf_{n \to \infty} \mle{\ell}(X_n)\geq \sqrt{2}\ell_1$ and $\liminf_{n \to \infty} \cle{\ell}(X_n)\geq \sqrt{2}\ell_1/4$.
The observational points within the interval $[-1,1]$ are equispaced, and we vary their number from 10 to 100. 
Due to its smoothness, the Gaussian kernel suffers from extreme numerical ill-conditioning; see ~\cite[Section~3.2]{KarvonenOates2026} and~\cite{barnett2024uniform}.
To resolve this, we employ high-precision arithmetic. Specifically, we use Julia's built-in \texttt{BigFloat} data type to perform calculations with up to 1408-bit precision. For optimisation, we utilise the \textsf{Particle Swarm} algorithm from \texttt{Optim.jl}, which supports \texttt{BigFloat}.

\begin{figure}[htbp]
    \centering
    \begin{tikzpicture}[scale = 0.95]
        \begin{groupplot}[
            group style={
                group size=1 by 3,
                vertical sep=1.2cm, 
                x descriptions at=edge bottom 
            },
            width=\linewidth,
            height=0.35\linewidth,
            xmin=10, xmax=100,
            grid=none,
            tick label style={font=\small \bfseries},
            label style={font=\small },
            title style={font=\bfseries, yshift=-0.8ex},
            xlabel={No. of Observational Points ($n$)},
            ylabel={ Estimates $\mathbf{\hat{\ell}^f}$ },
            legend style={font=\small, at={(0.98,0.95)}, anchor=north east}
        ]
            
            \nextgroupplot[
                title={$f(x) = k_{-2.0}^{1.7}(x)$},
                ymin=2.37, ymax=2.9,
                legend entries={ML Estimates, CV Estimates, ~$\ell^{*} = {\sqrt{2}\times1.7}$}
            ]
            \addplot[color=cyan, mark = *, mark size=0.7pt, line width=0.7pt] table{Data/mle1.dat};
            \addplot[color=orange,densely dashed, mark = *, mark size=0.7pt, line width=0.7pt] table{Data/cv1.dat};
            \addplot[color=blue,densely dashed, line width = 0.7pt, domain = 10:100]{2.40416305603};
            \nextgroupplot[
                title={$f(x) = k_0^{1.7}(x) +  k_{1}^{0.5}(x)$},
                ymin=0.7, ymax=0.815,
                legend entries={ML Estimates, CV Estimates, ~$\ell^{*}  = {\sqrt{2}\times0.5}$}
            ]
            \addplot[color=cyan, mark = *, mark size=0.7pt, line width=0.7pt] table{Data/mle2.dat};
            \addplot[color=orange, densely dashed, mark = *, mark size=0.7pt, line width=0.7pt] table{Data/cv2.dat};
            \addplot[color=blue,densely dashed, line width = 0.7pt, domain = 10:100]{0.70710678118};
            \nextgroupplot[
                title={$f(x) = k_0^{4.3}(x) + 2 k_{-1}^{5.0}(x) -3 k_{2}^{6.1}(x)$},
                ymin=6.0, ymax=7.05,
                legend entries={ML Estimates, CV Estimates, ~$\ell^{*} = \sqrt{2}\times4.3$}
            ]
            \addplot[color=cyan, mark = *, mark size=0.7pt, line width=0.7pt] table{Data/mle3.dat};
            \addplot[color=orange, densely dashed, mark = *, mark size=0.7pt, line width=0.7pt] table{Data/cv3.dat};
            \addplot[color=blue, densely dashed, line width = 0.7pt, domain = 10:100]{6.081118318204309};
        \end{groupplot}
    \end{tikzpicture}
    \caption{Lengthscale estimates as functions of the number of observation points $n$ for three response functions. The dashed line marks $\ell^{*} = \sqrt{2}\ell_1$, where $\ell_1$ is the smallest lengthscale used to define $f$ in~\eqref{eq:LinearSumofFeatures}.
    The scaled lengthscale $\ell^*$ determines the inclusion of $f$ in the RKHS of the Gaussian kernel, in that $f \in \mathcal{H}_{\ell_0}$ if $\ell_0 \in (0, \ell^*)$ and $f \notin \mathcal{H}_{\ell_0}$ if $\ell_0 \geq \ell^*$.
    The lengthscale estimates (both ML and CV estimates) are represented along the $y$-axis.}
    \label{fig:Sim_1}
\end{figure}
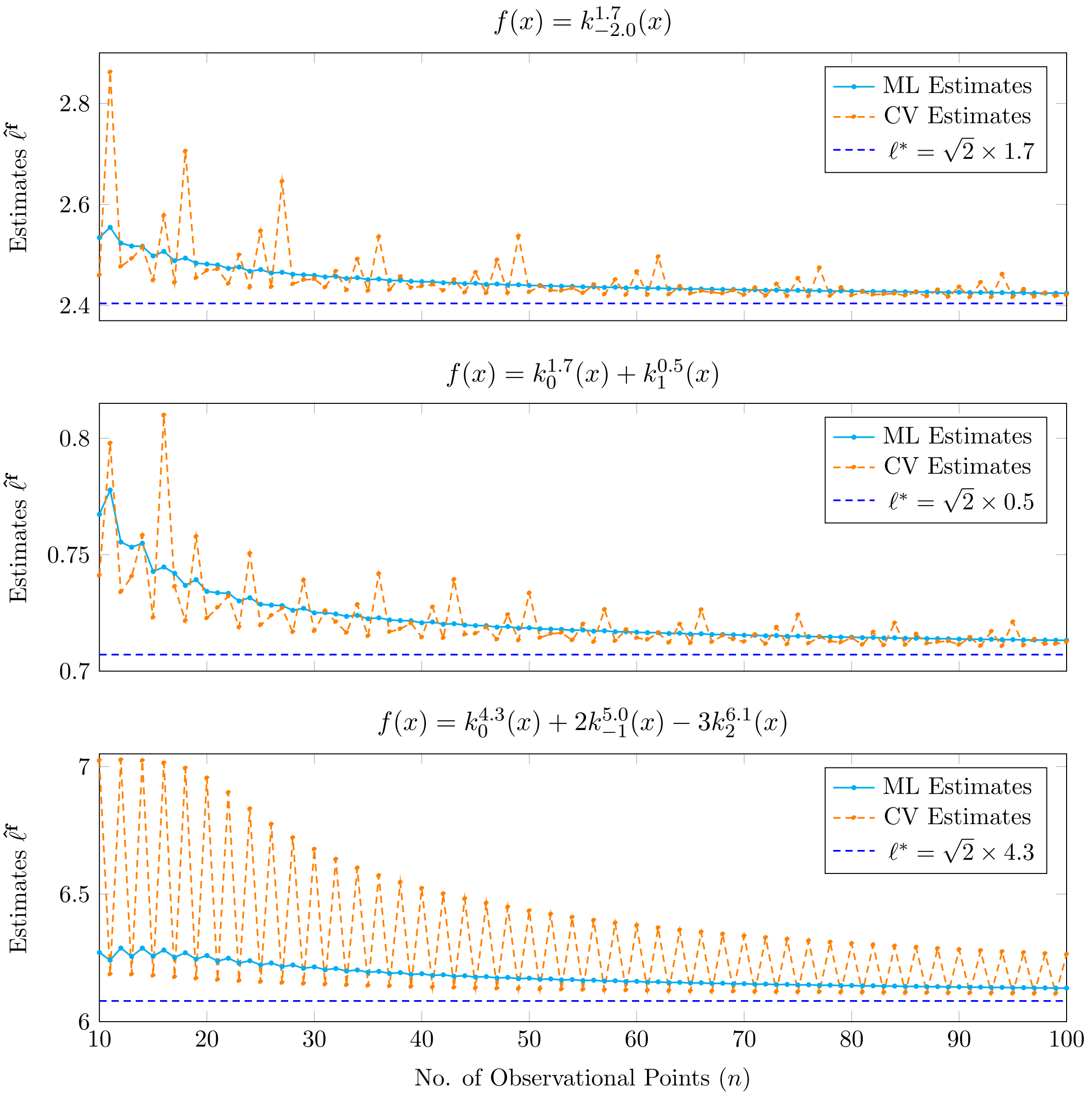

The asymptotic lower bounds for ML and CV estimates of the lengthscale given in Theorem~\ref{th:BdSubsetRKHS} are independent of the distribution of
the sampling points within $[-1, 1]$. While the results in Figure~\ref{fig:Sim_1} use equispaced points, we also ran the simulations using both Chebyshev nodes and a Van der Corput sequence. These configurations yielded results similar to those presented in Figure~\ref{fig:Sim_1} and are therefore not reported. Our numerical experiments indicate that both ML and CV estimators approach $\ell^{*} = \sqrt{2}\ell_1$. This suggests that the theoretical lower bound of CV estimator in Theorem~\ref{th:BdSubsetRKHS} may be suboptimal, as the estimator converges to the true parameter value. Overall, we find that numerical experiments of lengthscale estimation align with Theorem~\ref{th:BdSubsetRKHS}.

\section{Discussion}
This article proved lower bounds for the ML and CV estimators of lengthscale of a Gaussian process with a Gaussian kernel, analysing three distinct functional regimes as the number of data points $n$ approaches infinity. 
In the first case, when the target function is only assumed to belong to the RKHS of the Gaussian kernel for a given lengthscale, both estimators are bounded below by a constant multiple of that true value. 
In the second case the function is bandlimited and both estimators blow up at least as fast as $\sqrt{\log n}$. Finally, when the Fourier transform of the response function decays as $\exp(-a|\omega|^\gamma)$ for $\gamma > 2$ and $a>0$, the estimators blow up at least as fast as \smash{$(\sqrt{\log n } )^{{(\gamma-2)}/{\gamma}}$}.
The techniques we use are not refined enough to prove upper bounds on the estimators.
For a generic theory of upper bounds for Gaussian process parameter estimators, see Section 2.3.2 in~\cite{karvonen2023asymptotic}.
The general results established in Theorem~\ref{th:Gen_Method} can likely be applied to find lower bounds on lengthscale estimators for other kernels, such as Wendland, Gevrey~\cite[Section~4.3]{KarvonenOates2026}, and power series kernels~\cite{Zwicknagl2009}.
In this work, we focus on the Gaussian kernel on the interval $[-1, 1]$.
The tools we use---specifically the Gaussian Gram determinant bounds---only work in the univariate setting. It remains an open question how to generalise the results to higher dimensions. 

\section*{Acknowledgements}
The authors used Claude Opus 5.5 for proofreading. The authors assume responsibility for all content.


\bibliographystyle{siamplain}
\bibliography{references}

\end{document}